\documentclass [11pt] {article}

\usepackage {amssymb,latexsym}
\usepackage {amsmath}
\usepackage {citesort}

\title {BRST Cohomology and Phase Space Reduction\\
        in Deformation Quantization}

\author {{\bf
          Martin
          Bordemann\thanks{Martin.Bordemann@physik.uni-freiburg.de}~,
          \addtocounter{footnote}{2}
          Hans-Christian
          Herbig\thanks{Herbig@physik.uni-freiburg.de}~,
          Stefan
          Waldmann\thanks{Stefan.Waldmann@physik.uni-freiburg.de}
         } \\[3mm]
         Fakult\"at f\"ur Physik\\Universit\"at Freiburg \\
         Hermann-Herder-Str. 3 \\
         79104 Freiburg i.~Br., F.~R.~G \\[3mm]
        }

\date{FR-THEP-98/25 \\[1mm]
      Revised version 
      \\[1mm]
      September 1999
      }

\newcommand{\id} {{\mathsf {id}}}
\newcommand{\tr} {\mathop{\mathsf{tr}}}

\newcommand{\ad} {{\mathrm {ad}}}
\newcommand{\Ad} {\mathop{\mathrm {Ad}}}
\newcommand{\Hom} {{\mathsf {Hom}}}

\newcommand{\image} {\mathop{\mathrm{im}}}

\newcommand{\TinyW} {{\mbox{\rm \tiny W}}}
\newcommand{\weyl} {\mathop{\circ_\TinyW}}
\newcommand{\adw} {{{\mathrm {ad}}_\TinyW}}

\newcommand{\TinyS} {{\mbox{\rm \tiny S}}}
\newcommand{\std} {\mathop{\circ_\TinyS}}
\newcommand{\ads} {{{\mathrm {ad}}_\TinyS}}

\newcommand{\kcliff} {\mathop{\circ_\kappa}}
\newcommand{\adk} {{{\mathrm{ad}}_\kappa}}
\newcommand{\kstar}{\mathop{\star_\kappa}}

\newcommand{\im} {{\mathrm i}}
\newcommand{\eu} {{\mathrm e}}

\newcommand{\Grass} {\mbox{$\bigwedge$}}

\newcommand{\kbrsc} {\Theta_{\kappa}}
\newcommand{\wbrsc} {\Theta_{\mbox {\rm {\tiny{W}}}}}
\newcommand{\sbrsc} {\Theta_{\mbox {\rm {\tiny{S}}}}}

\newcommand{\kbrso} {\mbox{\boldmath $\mathcal D$}_{\kappa}}
\newcommand{\wbrso} {\mbox{\boldmath $\mathcal D$}_{\mbox {\rm {\tiny{W}}}}}
\newcommand{\sbrso} {\mbox{\boldmath $\mathcal D$}_{\mbox {\rm {\tiny{S}}}}}

\newcommand{\kpr} {\star_{\kappa}}
\newcommand{\wpr} {\star_{\mbox {\rm {\tiny{W}}}}}
\newcommand{\spr} {\star_{\mbox {\rm {\tiny{S}}}}}

\newcommand{\kos} {\partial}
\newcommand{\ce} {\delta}

\newcommand{\qce} {\mbox{\boldmath ${\delta}$}}
\newcommand{\qkos} {\mbox{\boldmath ${\partial}$}}

\newcommand{\qrestr} {\mbox{\boldmath $\iota^*$}}

\newcommand{\qmm} {\mbox{\boldmath $J$}}

\newcommand{\uni}{\mbox{\boldmath $u$}}

\newcommand{\quop}{\mbox{\boldmath $q$}}

\newcommand{\cuop}{\mbox{\boldmath $c$}}

\newcommand{\ceind} {{\mbox{\rm \tiny CE}}}
\newcommand{\kosind} {{\mbox{\rm \tiny Kos}}}
\newcommand{\brsind} {{\mbox{\rm \tiny BRST}}}

\newcommand{\smu} {\mbox{\boldmath $M$}_{\mbox{\tiny S}}}
\newcommand{\amu} {\mbox{\boldmath $M$}_{\mbox{\tiny A}}}
\newcommand{\wqce} {\mathop{\mbox{\boldmath $\delta$}_{\mbox {\rm \tiny W}}}}
\newcommand{\wqkos} {\mathop{\mbox{\boldmath $\partial$}_{\mbox{\rm \tiny W}}}}
\newcommand{\qhomot} {\mbox{\boldmath $h$}}
\newcommand{\qcohom} {\mbox{\boldmath $H$}}
\newcommand{\qideal} {\mbox{\boldmath $\mathcal I$}}
\newcommand{\qbideal} {\mbox{\boldmath $\mathcal B$}}
\newcommand{\QPsi} {\mbox{\boldmath $\Psi$}}
\newcommand{\qE} {\mbox{\boldmath $E$}}
\newcommand{\qW} {\mbox{\boldmath $W$}}

\newcommand{\ClassLieM} {{\varrho_{\mbox{\tiny $M$}}}}
\newcommand{\ClassLieC} {{\varrho_{\mbox{\tiny $C$}}}}
\newcommand{\QuantLieM} {{\mbox{\boldmath $\varrho$}_{\mbox{\tiny $M$}}}}
\newcommand{\QuantLieC} {{\mbox{\boldmath $\varrho$}_{\mbox{\tiny $C$}}}}

\newcommand{\starred} {\mathop{*_{\mbox{\rm \tiny red}}}}

\newcommand{\pervstar} {\mathop{\tilde *}}

\newcommand{\prol} {\mathop{\mathrm {prol}}}

\newenvironment {proof}{\small {\sc Proof:}}{{\hspace*{\fill} $\square$}}

\newtheorem {lemma} {Lemma}
\newtheorem {proposition} [lemma] {Proposition}
\newtheorem {theorem} [lemma] {Theorem}

\newtheorem {definition}[lemma] {Definition}
\newtheorem {remark}[lemma]{Remark}

\begin{document}

\maketitle
\thispagestyle{empty}

\begin{center}
\large
To the memory of \\
Mosh{\'e} Flato
\end{center}

\vspace{1cm}

\begin{abstract}
In this article we consider quantum phase space reduction 
when zero is a regular value of the momentum map.
By analogy with the classical case we define
the BRST cohomology in the framework of deformation quantization. We
compute the quantum BRST cohomology in terms of a `quantum'
Chevalley-Eilenberg cohomology of the Lie algebra on the constraint
surface. To prove this result, we construct an explicit
chain homotopy, both in the classical and quantum case, which is
constructed out of a prolongation of functions on the constraint
surface. We have observed the phenomenon
that the quantum BRST cohomology cannot always be used for quantum
reduction, because generally its zero part is no longer a
deformation of the space of all smooth functions on 
the reduced phase space. But in case the group action is `sufficiently
nice', e.g.~proper (which is the case for all compact Lie group
actions), it is shown for a strongly invariant star product that the
BRST procedure always induces a star product on the reduced phase
space in a rather explicit and natural way. Simple examples and
counter examples are discussed. 
\end{abstract}

\newpage
\thispagestyle{empty}
\tableofcontents

\newpage

%
%

%
%

\section {Introduction}
\label {IntroSec}

The aim of this article is to give a deformation quantization 
formulation of the method of BRST (Becchi-Rouet-Stora-Tyutin) cohomology
which has been introduced and is frequently used in physics for
the quantization of so-called constrained systems 
(see e.g.~\cite{HT92} and references therein): in symplectic geometry
these systems are known as reduced symplectic manifolds.

Deformation quantization has been introduced in
\cite{BFFLS78}. Quantization is formulated as an associative formal
deformation, the so-called star product, of the commutative
algebra of complex-valued $C^\infty$-functions on a symplectic
(or more generally, a Poisson) manifold by a formal series in 
$\lambda$ (which corresponds to Planck´s constant $\hbar$ in the
convergent case) of bidifferential operators such that the term of
order zero is pointwise 
multiplication and the commutator of the first order term is equal to
$\im$ times the Poisson bracket. By now there are general existence
and classification theorems for star products on every Poisson manifold 
(see \cite{Kon97b}, and earlier, for symplectic manifolds, 
\cite{DL83b,Fed94a,NT95a,BCG97}). Representation theory
for the deformed algebras in the spirit of $C^*$-algebras has been introduced
by \cite{BW98a} by formulating formally positive functionals and formal
GNS representations.

The reduction of symplectic manifolds by means of a sufficiently
`nice' Hamiltonian action of a Lie group $G$ (with Lie algebra 
$\mathfrak g$) has been formulated by Marsden and
Weinstein (see e.g. \cite[Chapter~4]{AM85}), and the method of BRST
cohomology has been transferred to symplectic geometry by 
Browning and McMullan \cite{McMul87,BMcM87}, 
Kostant and Sternberg \cite{KS87}, and Henneaux and Teitelboim
\cite{HT88}, see also \cite{FHST89,Dub87}.

Let us recall briefly the definition of those two objects:
\begin{itemize}
\item The starting point of reduction is an
      $\Ad^*$-equivariant momentum map $J: M \to \mathfrak g^*$ for
      the $G$-action (see e.g.~\cite[p.~276]{AM85}) whose level
      surface $C:=J^{-1}(\{0\})$ plays the r\^{o}le of the constraint
      surface in physics. 
      The reduced phase space $M_{\rm red}$ is then the symplectic manifold
      $C/G$ where usually the action of $G$ on $C$ is supposed to be
      proper and free in order to guarantee a compatible
      differentiable structure on $M_{\rm red}$. 
      In that case thanks to the Dirac method we can see
      the space $C^\infty (M_{\rm red})$ as a quotient 
      $\mathcal B_C \big/ \mathcal I_C$ where $\mathcal I_C$ is the
      vanishing ideal of $C$ and $\mathcal B_C$ is its
      normalizer with respect to the Poisson bracket.

\item For the BRST method the Poisson algebra $C^\infty(M)$ is
      tensored by the Gra{\ss}mann algebra over the direct sum of the
      Lie algebra and its dual (the latter space 
      $\Grass \mathfrak g^* \otimes_{\mathbb R} \Grass \mathfrak g$
      itself becomes a super Poisson algebra by means of the natural
      pairing between $\mathfrak g$ and $\mathfrak g^*$) to an
      extended super Poisson algebra $\mathcal A$ which is called the
      classical BRST algebra. Roughly speaking, the super Poisson
      bracket of a suitable linear combination of the Lie structure
      and the momentum map, $\Theta$, serves as an odd Hamiltonian
      super-derivation of square zero, the so-called classical BRST
      operator, in the extended super Poisson algebra 
      $\mathcal A$. There is in addition the so-called ghost number or
      total degree derivation on $\mathcal A$ which is equal to $k-l$
      on each subspace 
      $\Grass^k \mathfrak g^* \wedge \Grass^l \mathfrak g \otimes C^\infty(M)$. 
      It turns out that the classical BRST operator is the total
      differential of a double complex whose vertical differential is
      a Chevalley-Eilenberg differential whereas the horizontal
      differential is twice the classical Koszul boundary operator.

      BRST cohomology gives then a very nice method to
      describe the space $C^\infty (M_{\rm red})$, a method that we
      shall use for quantization.
      In case $0$ is a regular value of the momentum map
      (and the action of $G$ on $M$ is nice: for example
      proper and free) the `ghost number zero part' of the
      cohomology of the classical BRST 
      operator is known to be isomorphic (as a Poisson algebra) to the
      space $C^\infty(M_{\rm red})$.
\end{itemize}

Recently there have been several articles in which phase space reduction
has been dealt with in deformation quantization: Fedosov has formulated
symplectic reduction in his scheme of star products for $U(1)$-actions
\cite{Fed94b} and for general compact groups \cite{Fed98a}. In particular
situations, phase space reduction methods have been used to compute 
explicit formulas for star products on all complex projective spaces
\cite{BBEW96a,BBEW96b}, on Gra{\ss}mann manifolds
\cite{Schir97,Schir98}, and for one-dimensional Lie algebras
\cite{Gloe98a,Gloe98b}. The method of BRST cohomology has
been successfully formulated for geometric 
quantization (see e.g.~\cite{Kim92} and \cite{Tuy92}, and references
therein), but there does not yet seem to be a treatment of BRST cohomology
in deformation quantization although BRST cohomology ``seems to be
well-suited to the recent work on the deformation approach to
quantization'', see \cite[p.~428]{McMul87}.

In this article, we shall give a quantum version of the
BRST method described before to get similar results for deformation
quantization on constrained systems. But we shall not restrict our
study only to the nice cases (such as regular value of the momentum
map, proper and free group action). Actually, we are convinced that a
treatment of BRST cohomology in deformation quantization has several
advantages:

Firstly, physicists using BRST
cohomology methods often complain about operator ordering problems
which forces them to {\em a priori assume the existence} of a quantum
BRST (cohomology) operator: to quote Henneaux and Teitelboim 
\cite[p.~297]{HT92}: 
``It will be assumed that one can find a charge $\Omega$ satisfying the 
nilpotence and hermiticity conditions [\ldots] Unlike the situation in the 
classical case, there is no {\em a priori} guarantee that this can always be 
done starting from a classical theory for which $[\Omega,\Omega]=0$, since 
the ordering of the factors comes in crucially.'' 
In contrast to that, deformation quantization
can also be viewed as a theory consistently overcoming and even classifying 
{\em all} possible operator orderings in situations where differential 
operator representations of the deformed algebra (for example in a
symbol calculus on cotangent bundles, see
e.g. \cite{BNW98a,BNW99a,Pfl98b,BNPW98} for a treatment on curved
configuration spaces) are possible. Moreover, there are general theorems 
in deformation quantization about the quantization of proper Hamiltonian
group actions \cite[p.~180--183]{Fed96}.

Secondly, in deformation quantization (as in the $C^*$-algebra theory)
the observable algebra is the principal object whereas representations are
subordinate. Therefore it is rather natural to check whether a classical 
BRST operator simply remains a cohomology operator when the super Poisson
bracket is replaced by the super-commutator of a $\mathbb Z_2$-graded
star product.

Thirdly, by its very definition deformation quantization allows
us to control the classical limit after the quantum reduction which often
ends up (in other quantization schemes) with abstract quotient algebras.

In this article we have come to the following principal results:
\begin{enumerate}
\item For every star product $*$ on $M$ covariant under the group
      action (which can be achieved for every proper Lie group action), i.e.
       \[  
           \langle J ,\xi\rangle * \langle J, \eta\rangle - 
           \langle J, \eta\rangle * \langle J, \xi\rangle 
           = \im\lambda \langle J, [\xi,\eta]\rangle 
           \qquad
           \forall \xi,\eta\in\mathfrak g,
      \]
      where $\langle \cdot, \cdot \rangle$ denotes the natural pairing,     
      we construct a one-parameter family of
      formal associative deformations of the classical extended super
      Poisson algebra, 
      $({\mathcal A}[[\lambda]],\kpr)_{\kappa \in [0,1]}$,
      which are all 
      equivalent by an explicit equivalence transformation $S_\kappa$ 
      and for which the corresponding BRST charge $\Theta_\kappa$ has square
      zero (Section \ref{QuantBRSTSec}).

\item We compute the family of quantum BRST operators $\kbrso$ in 
      $({\mathcal A}[[\lambda]],\kpr)$, i.e.~the graded star product
      commutators with $\Theta_\kappa$. It turns out that the quantum
      BRST operator which seems to be the most `natural' from the
      point of view of Clifford algebras (which we called
      the Weyl-ordered BRST operator $\wbrso$, and which
      corresponds to $\kappa=1/2$) 
      looks rather complicated and not very encouraging concerning
      cohomology computations. But luckily the quantum BRST operator
      corresponding to $\kappa=0$ (which we called standard ordered
      BRST operator $\sbrso$) ---which is conjugate by $S_{1/2}$ to
      $\wbrso$--- again defines a double complex 
      (Section~\ref{QuantBRSTSec}) thus 
      giving rise to deformed versions of the classical Koszul
      boundary operator and the classical Chevalley-Eilenberg
      differential.

\item For every regular level-zero surface $C$ 
      we compute the quantum BRST cohomology of
      $({\mathcal A}[[\lambda]],\sbrso)$ (which is again a 
      $\mathbb Z$-graded  associative algebra) by means of a 
      \emph{deformed augmentation}, i.e. a quantised version of the
      linear map restricting functions on the manifold to $C$: the
      result is that the quantum BRST cohomology is isomorphic in an
      explicit way (using deformed versions of classical chain
      homotopies of the classical Koszul complex) to the
      Chevalley-Eilenberg cohomology of the Lie algebra $\mathfrak g$
      with $\mathfrak g$-module $C^\infty(C)[[\lambda]]$ where the
      representation is a deformation of the usual Lie derivative of
      the vector fields of the classical 
      $\mathfrak g$-action. Moreover the quantum BRST 
      cohomologies of the above quantum BRST operators $\kbrso$ are
      all isomorphic as associative algebras to the cohomology of the
      standard ordered operator. Finally, we also arrive at the
      isomorphy of the quantum BRST cohomology algebra and a
      Dirac-type picture: we define a deformed version of the
      vanishing ideal $\mathcal I_C$ of the constraint surface, a
      certain left ideal $\qideal_C$ of 
      $(C^\infty (M)[[\lambda]], *)$, 
      and its idealiser $\qbideal_C$ modulo $\qideal_C$
      turns out to be naturally isomorphic to the quantum BRST 
      algebra (Section \ref{BRSTCohoSec}).

\item The natural question arising in view of the preceding result is
      the following: even in the nice case (regularity,
      proper and free action), does this deformed
      Chevalley-Eilenberg cohomology for ghost number zero reflect in
      an isomorphic manner the space of functions on the reduced phase
      space (in case this space exists)? As we shall show in Section
      \ref{ExampleSec} by a simple example dealing with a Hamiltonian
      circle action the answer is in general \emph{no}! It may happen
      that momentum map and star product are so `ill-adjusted' that
      the zeroth quantised Chevalley-Eilenberg cohomology of
      $\mathfrak g$ on the constraint surface $C$ becomes `much
      smaller' than the classical cohomology which is in bijection
      with $C^\infty(M_{\rm red})$.

\item However, as we show in Section \ref{ProGroupSec} there are large
      classes of examples in which the afore-mentioned pathology does
      not occur: the first class is the family of proper Hamiltonian
      $G$-spaces for which there always exist strongly invariant
      star products and $G$-equivariant chain homotopies and
      prolongations. Here, using such a strongly invariant
      star product, classical and quantum Chevalley-Eilenberg
      cohomologies of the Lie algebra $\mathfrak g$ on the constraint
      surface are simply equal. Moreover we get fairly explicit
      formulas for the star product on the reduced phase space in
      terms of the star product on the original symplectic manifold,
      an equivariant prolongation map and a deformed restriction
      map. This formula is particularly simple for global
      $G$-invariant functions thus serving to quantize integrable
      systems obtained by reduction.  
      The second class of examples consists of those situations where
      the first classical Chevalley-Eilenberg cohomology group of the
      Lie algebra $\mathfrak g$ on the constraint surface is zero. For 
      particular cases this is satisfied when the first de Rham
      cohomology of the Lie group vanishes. 
\end{enumerate}
The paper is organised as follows: in the first three Sections 
\ref{StarProdSec}, \ref{GeomClassSec}, and \ref{ClassBRSTSec} we
recall basic concepts and results in deformation quantization,
geometry and Koszul complex for constraint surfaces, and classical
BRST theory, respectively. 
As stated above, the main results of this paper are contained in Section 
\ref{QuantBRSTSec}, Section \ref{BRSTCohoSec}, Section \ref{ExampleSec}, and 
Section \ref{ProGroupSec}. In Section \ref{OutSec} we give a short
conclusion and discuss further problems and questions arising with our 
approach.

\vspace{0.5cm}

\noindent
Notation: Tensor products $\otimes$ are usually taken over 
$\mathbb C$. Otherwise the ring will be indicated as subscript. 
Moreover, $C^\infty (M)$ always denotes the space of smooth
complex-valued functions on $M$. Finally, $H^\bullet$ indicates a
$\mathbb Z$-grading of a module 
$H$ and analogously $\mathcal A^{\bullet,\bullet}$ denotes a 
$\mathbb Z \times \mathbb Z$ grading. A homogeneous map $\Phi$ of
degree $k$ is denoted by $\Phi: H^\bullet \to H^{\bullet+k}$.

%
%

\section {Star products and Hamiltonian Lie group/algebra actions}
\label {StarProdSec}

In this section we shall recall some basic concepts of deformation
quantization and Hamiltonian Lie group and Lie algebra actions in order to
establish our notation, see also e.g.~\cite{AM85}.

We consider a Poisson manifold $(M, \Lambda)$, i.e.~a smooth manifold $M$
with a Poisson tensor field $\Lambda \in \Gamma^\infty (\bigwedge^2 TM)$
such that the Schouten bracket $[\Lambda, \Lambda]$ vanishes, see
e.g.~\cite{CW99}. Then
$\{f, g\} := \Lambda(df, dg)$ defines a Poisson bracket which turns
$C^\infty (M)$ into a Poisson algebra. Here we always consider
\emph{complex-valued} functions and tensor fields. The vector field
$X_f := \Lambda (df, \cdot)$ is called the Hamiltonian vector field of
$f \in C^\infty (M)$. A particular case of a Poisson manifold is a
symplectic manifold $(M, \omega)$ where the symplectic form 
$\omega \in \Gamma^\infty (\bigwedge^2 T^*M)$ is a closed,
non-degenerate two-form. In this case the Hamiltonian vector field of
$f$ is defined by $i_{X_f} \omega = df$ and the Poisson bracket is
$\{f, g\} = \omega (X_f, X_g)$ whence the Poisson tensor $\Lambda$ is
just the `inverse' of $-\omega$.

Now we consider the space $C^\infty (M)[[\lambda]]$ of formal power series
in a formal parameter $\lambda$ as $\mathbb C[[\lambda]]$-module. Then a
\emph{star product} $*$ for $(M, \Lambda)$ is a 
$\mathbb C[[\lambda]]$-bilinear,
associative deformation of the pointwise product of $C^\infty (M)$ such that
\begin {equation} \label {StarProduct}
    f * g = \sum_{r=0}^\infty \lambda^r C_r (f, g),
\end {equation}
where $C_0 (f, g) = fg$, $C_1 (f, g) - C_1 (g, f) = \im \{f, g\}$, and all
$C_r$ are bidifferential operators vanishing on constants whence
$1*f = f = f*1$, see \cite{BFFLS78}. Sometimes further requirements
are made by specifying
certain parity or reality condition for the $C_r$. Furthermore a star
product is of \emph{Vey type} if the bidifferential operator $C_r$ is
of order $r$ in each argument for all $r$. One might also take 
\emph{local} operators $C_r$ instead of bidifferential ones but we
shall deal only with bidifferential ones for simplicity.
The formal parameter $\lambda$ plays the r{\^o}le of Planck's constant
$\hbar$ and may be substituted by $\hbar$ in convergent situations. 
The existence of such star
products was shown in the symplectic case by DeWilde and Lecomte \cite{DL83b},
Fedosov \cite{Fed86,Fed94a}, and Omori, Maeda, and Yoshioka \cite{OMY91}, and
recently by Kontsevich in the general Poisson case \cite{Kon97b}.
Two star products $*$ and $*'$ are called \emph{equivalent} if there exists a
formal series of differential operators
$T = \id + \sum_{r=1}^\infty \lambda^r T_r$ such that $T (f*g) = Tf *' Tg$
for all $f, g \in C^\infty (M)[[\lambda]]$.
The classification up to equivalence was done by Nest and Tsygan
\cite{NT95a,NT95b}, Deligne \cite{Del95},
Bertelson, Cahen and Gutt \cite{BCG97}, and Kontsevich \cite{Kon97b}.

Let $\mathfrak g$ be a finite-dimensional real Lie algebra with dual space
$\mathfrak g^*$. Recall that a \emph{Lie algebra action} of $\mathfrak
g$ on $M$ is a linear anti-homomorphism $\xi \mapsto \xi_M$ of
$\mathfrak g$ into the Lie algebra of all vector fields on $M$. It
follows that these vector fields define a representation $\ClassLieM$ 
of $\mathfrak g$ in the space $C^\infty(M)$ by
\begin{equation}
\label{KlassLieM}
  \ClassLieM(\xi)(f) := -\xi_M(f).
\end{equation}
Let $G$ be a Lie group with Lie algebra $\mathfrak g$. Any (left)
\emph{Lie group action} of 
$\Phi: G\times M \to M: (g,x) \mapsto \Phi_g(x)$ defines a Lie
algebra action by means of its infinitesimal generators
$\xi_M:=\left.\frac{d}{dt} \Phi_{\exp (t\xi)}\right|_{t=0}$. 
Recall that a Lie algebra action (or a Lie group action) on a Poisson
manifold $(M,\Lambda)$ is called Hamiltonian (see
e.g.~\cite[Sect.~4]{AM85} for details) if and only if there is a 
\emph{momentum map} of the action, i.e.~a $C^\infty$-map 
$J: M \to \mathfrak g^*$ such that for every $\xi \in \mathfrak g$ 
\begin{equation} 
\label{MomentMap}
    X_{\langle J,\xi\rangle} = \xi_M.
\end{equation}
Moreover we require equivariance of $J$ with respect to the coadjoint
representation of $\mathfrak g$ and $G$, i.e.
$\{\langle J,\xi \rangle, \langle J,\eta \rangle\}=\langle J,[\xi,\eta] \rangle$
for all $\xi,\eta\in\mathfrak g$ in the case of a Hamiltonian Lie algebra
action and $J(\Phi_g (x)) = \Ad^*(g) J(x)$ for all $g \in G$
and $x \in M$ in case of a Hamiltonian Lie group action. In the last case the
quadruple $(M, \Lambda, G, J)$ is usually called a \emph{Hamiltonian $G$-space}.
We shall speak of a \emph{Hamiltonian $\mathfrak g$-space} 
$(M,\Lambda,\mathfrak g, J)$ in the more general case of a Hamiltonian
Lie algebra action. In the symplectic case we shall denote this by
$(M, \omega, \mathfrak g, J)$ and $(M, \omega, G, J)$, respectively.

For a physically reasonable quantization procedure one certainly has to
impose more conditions on a star product beside the defining ones, since
particular properties of the Poisson manifold, as e.g.~symmetries, should
be preserved under quantization. This leads to the various definitions of
`invariance' for star products under a given classical Lie group or Lie
algebra action: In the context of deformation quantization the
following notions of invariance are commonly used, see 
e.g.~\cite{ACMP83}: The star product $*$ is called 
\begin {itemize}
\item \emph{invariant} if 
      $\Phi_g^* (f*h) = \Phi^*_g (f) * \Phi^*_g (h)$ for a Hamiltonian Lie group
      action and, more generally, for a Hamiltonian Lie algebra action
      $\{ \langle J, \xi\rangle, f*h\} = 
      \{ \langle J, \xi\rangle, f\} * h + 
      f * \{\langle J, \xi\rangle, h\}$ for all
      $g \in G$ resp. $\xi \in \mathfrak g$ and
      $f, h \in C^\infty (M)[[\lambda]]$,
\item \emph{covariant} if 
      $\langle J, \xi \rangle * \langle J, \eta\rangle - 
      \langle J, \eta \rangle * \langle J, \xi \rangle = 
      \im \lambda \langle J, [\xi, \eta]\rangle$
      for all $\xi, \eta \in \mathfrak g$ for both types of Hamiltonian
      action, and finally,
\item \emph{strongly invariant} if 
      $\langle J,\xi \rangle * f - f * \langle J, \xi\rangle 
      = \im \lambda \{\langle J, \xi \rangle, f\}$
      for all $\xi \in \mathfrak g$ and $f \in C^\infty (M)[[\lambda]]$
      for both types of Hamiltonian action.
\end {itemize}
Then clearly strong invariance implies both invariance and covariance
in case the Lie group $G$ is connected. Furthermore one can allow
quantum corrections to the momentum map leading to the notion of a
\emph{quantum momentum map}, see e.g.~\cite{Xu98,BBEW96a}. We 
consider a formal series
$\qmm = \sum_{r=0}^\infty \lambda^r \qmm_r: M \to \mathfrak g^*[[\lambda]]$
of smooth functions $\qmm_r : M \to \mathfrak g^*$ such that
$\qmm_0 = J$ is the classical momentum map and $\qmm$ satisfies
$\langle \qmm, \xi \rangle * \langle \qmm, \eta\rangle 
- \langle \qmm, \eta \rangle * \langle \qmm, \xi\rangle =
\im\lambda \langle \qmm, [\xi, \eta]\rangle$
for all $\xi, \eta \in \mathfrak g$. In this case the star product will
be called \emph{quantum covariant}, and clearly covariance 
with respect to $J$ implies quantum covariance 
for $\qmm = J$. Moreover quantum covariance implies 
that the Lie algebra $\mathfrak g$ acts by `inner' derivations on the
algebra $(C^\infty (M)[[\lambda]], *)$, where the representation
$\QuantLieM$ is given by
\begin{equation}
\label{QuantLieMDef}
    \QuantLieM (\xi) = 
    \frac{1}{\im\lambda} \ad_* (\langle \qmm, \xi\rangle) 
\end{equation}
for $\xi \in \mathfrak g$. Here $\ad_*$ stands for taking commutators
with respect to the star product $*$. This motivates the following definition:
\begin{definition}
\label{QuantumGSpace}
The quadruple $(M, *, \mathfrak g, \qmm)$ is called Hamiltonian
quantum $\mathfrak g$-space if $(M, *)$ is a Poisson manifold with
star product such that $\qmm$ is a quantum momentum map and $*$ is
quantum covariant under $\mathfrak g$.
\end{definition}
Finally a quantum $\mathfrak g$-space 
$(M, *, \mathfrak g, \QuantLieM)$ is defined to be a Poisson manifold
with star product such that the Lie algebra $\mathfrak g$ acts via
$\QuantLieM$ by not necessarily inner star product derivations on 
$C^\infty (M)[[\lambda]]$. Given a Hamiltonian quantum $\mathfrak g$-space 
$(M, *, \mathfrak g, \qmm)$ we call 
$(M, \Lambda, \mathfrak g, J = \qmm_0)$ the corresponding 
\emph{classical limit}.

We shall now consider particular Hamiltonian
group actions which all imply the existence of a strongly invariant
star product in the 
symplectic case. Recall that a Lie group action $\Phi:G\times M\rightarrow M$
is called \emph{proper} if the map 
$\hat{\Phi}:G\times M \to M \times M: (g,m) \mapsto (\Phi(g,m),m)$ is
proper, i.e. inverse images of compact sets are compact. 
Equivalently, for any sequences $(x_n)_{n\in\mathbb N}$ in $M$ and
$(g_n)_{n\in\mathbb N}$ in $G$ such that $x_n\to x$ and
$\Phi(g_n,x_n)\to y$ there is a subsequence of $(g_n)_{n\in\mathbb N}$
converging to $g\in G$ such that $y=\Phi(g,x)$. Proper group actions
always have closed orbits and compact isotropy groups (see also
\cite[Chap.~4]{AM85}). One has the following implications in case 
$(M, \omega)$ is symplectic: 
\begin {equation} \label {VariousGroupActions}
\begin {minipage}{1.5cm}
\begin {center}
compact group action
\end {center}
\end {minipage}
\hfill
\Longrightarrow
\hfill
\begin {minipage}{1.5cm}
\begin {center}
proper group action
\end {center}
\end {minipage}
\hfill
\Longrightarrow
\hfill
\begin {minipage}{3cm}
\begin {center}
connection preserving group action
\end {center}
\end {minipage}
\hfill
\Longrightarrow
\hfill
\begin {minipage}{3cm}
\begin{center}
existence of strongly invariant star product
\end {center}
\end {minipage}
\end {equation}
While the first two implications are well-known general geometrical
results on Lie group actions, see e.g.~\cite [Thm.~4.3.1] {Pal61}, the
existence of a strongly invariant star product in the case where the
group leaves invariant a connection is proved by Fedosov's techniques
\cite [Sect.~5.8] {Fed96}. 

%
%

\section {Geometry of constraint surfaces and classical
          phase space reduction}
\label {GeomClassSec}

In this section we shall briefly recall the relation between
the geometry of constraint surfaces in a (Poisson) manifold
$M$ and certain subspaces of smooth functions on $M$.

Let $C$ be a regular closed sub-manifold of a manifold $M$
and denote by $\iota: C \to M$ the canonical embedding. Let
$\mathcal I_C$ denote the \emph{vanishing ideal of $C$},
i.e.~the subspace of $C^\infty (M)$ of all those functions
which vanish on $C$. The following well-known lemma shows the
existence of a \emph{prolongation} of smooth complex-valued functions
on $C$ in an open neighbourhood of $C$:
\begin{lemma}
\label{L31}
Let $C$ be a regular closed sub-manifold of a manifold $M$.
Then there is an open neighbourhood $U$ of  $C$ and a subspace
$\mathcal F_C$ of $C^\infty(M)$ such that:
\begin{enumerate}
\item Each $f \in \mathcal F_C$ is supported in $U$.
\item The restriction $\iota^*: \mathcal F_C \to C^\infty (C)$
      is a bijection. 
      We shall call its inverse $\prol$ the prolongation of
      $\phi \in C^\infty(C)$ to $M$. 
 \item The space $C^\infty(M)$ decomposes into the direct sum
      $\mathcal F_C \oplus \mathcal I_C$.
\item Let $\phi\in C^\infty(C)$ have compact support. Then
      $\prol\phi$ has compact support in $M$.
\end{enumerate}
In particular this entails that the map $\iota^*$ induces a
canonical bijection of the quotient $C^\infty(M)/\mathcal
I_C$ onto $C^\infty(C)$.
\end{lemma}
Note that we have borrowed the notation
$\prol$ from a paper by Gl\"{o}{\ss}ner \cite{Gloe98a,Gloe98b}, but
Gl\"{o}{\ss}ner denotes the projection $\prol\iota^*$ by $\prol$.

Suppose next that there is a proper action of a Lie group
$G$ on $M$ such that the sub-manifold $C$ is invariant under
this action. Then we have the following
\begin{lemma}
\label{L31inv}
Suppose that a Lie group $G$ properly acts on the manifold
$M$ such that the sub-manifold $C$ is preserved under this
action. Then the open neighbourhood $U$ and the subspace
$\mathcal F_C$ satisfying the properties of the preceding
Lemma can in addition be chosen to be invariant under this
action such that the prolongation intertwines the action on
$C$ with the action on $M$.
\end{lemma}
\begin{proof}
The existence of a $G$-invariant prolongation map is shown by using
a \emph{$G$-invariant tubular neighbourhood} $U$ of $C$ in $M$: consider
the conormal bundle $E:=\{\alpha\in T^*M|_C \;|\; \alpha(v)=0~\forall v\in 
T_{\tau(\alpha)}C \}$ (where $\tau$ denotes the cotangent bundle 
projection) on which $G$ acts in a canonical manner such that $\tau$
is $G$-equivariant; a $G$-invariant tubular neighbourhood consists in
the following data: a $G$-invariant open neighbourhood $N$ of the zero
section $C$ in $E$ and a $G$-equivariant diffeomorphism $\Phi$ of $N$ onto
a $G$-invariant open neighbourhood of $C$ in $M$ restricting to the identity
map on the zero section. We denote the $G$-equivariant pushed-forward
projection $\tau\circ \Phi^{-1}$ by $\tilde{\tau}$. We shall postpone
a sketch of an existence proof of the $G$-invariant tubular neighbourhood
at the end of this proof.

Consider now a $G$-invariant smooth partition of unity $\psi_U+\psi_W=1$ subordinate
to the $G$-invariant open covering of $M$ by $U$ and $W:=M\setminus C$
(which exists thanks to the properness of the $G$-action, see
\cite[p.~78, Thm.~5.2.5.]{PT88}).
For any $\phi\in C^\infty(C)$ define $\prol \phi$ to be equal to
$\psi_U (\phi\circ \tilde{\tau})$ on $U$ and zero outside of $U$. Clearly,
$\prol$ is $G$-equivariant and satisfies the asserted properties.

The existence proof of a $G$-invariant tubular neighbourhood can largely
be copied from the case $G=\{e\}$ in Lang's book \cite[p.~108--110]{Lang95}
by observing the following additional facts: 

Thanks to the properness of the $G$-action there is a $G$-invariant Riemannian 
metric on $M$ (see \cite[p.~316, Thm~4.3.1.]{Pal61}) inducing a $G$-equivariant
vector bundle isomorphism of $E$ onto the Riemannian normal bundle of
$C$ in $TM|_C$. $\Phi$ can then be defined as this morphism followed by the 
exponential map
of the metric, which can easily be seen to be a well-defined $G$-equivariant
local diffeomorphism of a $G$-invariant open neighbourhood of $C$ in $E$
onto a $G$-invariant open neighbourhood of $C$ in $M$.

It is more difficult to make $\Phi$ injective on a possibly smaller $G$-invariant
open neighbourhood of $C$ in $E$: suppose that for each $c\in C$ there is a 
\emph{$G$-invariant open neighbourhood}
restricted to which $\Phi$ is injective (a fact which we shall show further down).
Then the method of patching together local inverses (which are necessarily
$G$-equivariant) explained in Lang's book (following an argument by Godement)
can be transferred to our case since by properness of the $G$-action the quotient space 
$M/G$ is still paracompact (although in general no longer a manifold), see 
\cite[p.~302, Prop.~1.2.8. and p.~316, Thm~4.3.1.]{Pal61}: this guarantees
the existence of $G$-invariant locally finite coverings allowing for
$G$-invariant shrinkings and closures making all sets appearing in the
standard tubular neighbourhood proof $G$-invariant.

Finally, suppose that there were a point $c\in C$ having no $G$-invariant
open neighbourhood in $E$ restricted to which $\Phi$ is injective.
Then there would be a sequence $(X_n)$ of open neighbourhoods of $c$ having
compact closure and intersection $\{c\}$, two sequences 
$a_n,b_n\in X_n\setminus \{c\}$ converging both to $c$, and a sequence $g_n\in G$ 
such that $g_n a_n \neq b_n$ but $g_n\Phi(a_n)=\Phi(g_n a_n) = \Phi(b_n)$ for
all positive integers $n$ (using the open $G$-invariant sets $\{gx|g\in G, x\in X_n\}$).
Since $G$
properly acts on $M$ we can assume (by possibly restricting to a subsequence)
that $g_n$ converges to $g\in G$ with $gc=c$, hence $g_n b_n$ converges to $c$.
But then there is an integer $n_0$ such that $a_{n_0}$ and $g_{n_0}b_{n_0}$ both are
in a (not necessarily $G$-invariant) tubular neighbourhood $N$ of $C$ in $E$
restricted to which $\Phi$ is injective which is a contradiction.
\end{proof}

\noindent A prolongation constructed in the above proof will be called a
\emph{geometric} prolongation.

Let us suppose from now on that the manifold $M$ is
\emph{symplectic} with symplectic form $\omega$, and that 
$(M,\omega,\mathfrak g,J)$ is a Hamiltonian $\mathfrak g$-space. In
the rest of this paper we shall very often encounter the situation
that $0$ is a regular value of the equivariant momentum map $J$ and
that the \emph{constraint surface} $C := J^{-1}(\{0\})$ is nonempty in
which case we shall henceforth call the quintuple 
$(M, \omega, \mathfrak g, J, C)$ a \emph{Hamiltonian 
$\mathfrak g$-space with regular constraint surface} which in physics
is often called the irreducible case. 
The vanishing ideal $\mathcal I_C$ is a Poisson sub-algebra of 
$C^\infty(M)$, which is equivalent to $C$ being a coisotropic
(`first class') sub-manifold. 
Let $\mathcal B_C$ denote the \emph{normalizer} of the
vanishing ideal, i.e.~the space of all those functions
in $C^\infty(M)$ whose Poisson bracket with every
function in the vanishing ideal is again contained in
the vanishing ideal. Then $\mathcal I_C$ is a Poisson
ideal in $\mathcal B_C$ and the quotient 
$\mathcal B_C/\mathcal I_C$ becomes a Poisson algebra, see 
e.g.~\cite[p.~443]{Kim92}, \cite[p.~417--418]{AM85} or
Gl{\"o}{\ss}ner's paper \cite{Gloe98a} for a proof.
For Hamiltonian $G$-spaces we
shall speak of \emph{Hamiltonian $G$-spaces with regular constraint
surface} if the corresponding Hamiltonian $\mathfrak g$-space has a 
regular constraint surface. In this case it is known that
the group action on $C$ has zero-dimensional isotropy groups.
In order to define a smooth manifold structure on the reduced
phase space $M_{\rm red} := C/G$ which will turn the canonical
projection $\pi: C \to M_{\rm red}$ into a smooth open submersion, 
the group action does not have to be proper on all of $M$ but has to
be `sufficiently nice', e.g. proper and free, on the constraint
surface $C$ only. The following description of the space of
smooth complex-valued functions on the reduced space is
well-known:
\begin{proposition}
\label{ReducClas}
Let $(M, \omega, \mathfrak g, J, C)$ be a Hamiltonian $G$-space with
regular constraint surface, such that the connected Lie group $G$ acts 
in a sufficiently nice way (e.g. freely and properly) on $C$ such that 
$M_{\rm red}$ exists and the canonical projection $\pi$ is a smooth
open submersion.
\begin{enumerate}
\item The maps $\pi^*: C^\infty(M_{\rm red}) \to C^\infty (C)$
      and $\iota^*: C^\infty (M) \to C^\infty (C)$ induce the
      following bijections on the space of all $G$-invariant functions 
      $C^\infty (C)^G$
      \begin{equation}
      \label{NeueGleichung}
          \pi^*: C^\infty (M_{\rm red}) 
          \stackrel{\simeq}{\longrightarrow} C^\infty (C)^G,
          \qquad
          \iota^*: \mathcal B_C \big/ \mathcal I_C
          \stackrel{\simeq}{\longrightarrow} C^\infty (C)^G.
      \end{equation}

\item For any chosen prolongation $\prol$ one has      
      $\prol\pi^*C^\infty(M_{\rm red})\subset \mathcal B_C$.
      Moreover, the space of all $G$-invariant functions 
      $C^\infty(M)^G$ on $M$ is contained in $\mathcal B_C$.

      In case the group action is proper on all of $M$ we have that
      $\iota^*C^\infty(M)^G = \iota^*\mathcal B_C$.

\item Suppose that the Lie group $G$ acts properly and freely on $C$
      Then there is a left inverse 
      $\sigma: C^\infty (C) \to C^\infty (M_{\rm red})$ of $\pi^*$,
      i.e.~$\sigma\pi^* = \id_{C^\infty (M_{\rm red})}$ (this is the
      `gauge fixing map').
\item The Poisson bracket of two functions 
      $\phi_1, \phi_2 \in C^\infty (M_{\rm red})$ can be written as
      \begin{equation}
      \label{PoissonUnten}
          \{\phi_1, \phi_2\}_{\rm red} =
          \sigma \iota^* \{\prol\pi^*\phi_1, \prol\pi^*\phi_2\}
          \quad
          \textrm{or}
          \quad
          \pi^*\{\phi_1, \phi_2\}_{\rm red} = 
          \iota^* \{\prol\pi^*\phi_1, \prol\pi^*\phi_2\}.
      \end{equation}     
\end{enumerate}
\end{proposition}
\begin{proof}
The first point is a direct consequence of the definitions and of 
Lemma \ref{L31}. To prove that any function in $C^\infty (M)^G$ lies
in $\mathcal B_C$ we can use an explicit description of $\mathcal I_C$ 
generated by $\langle J,\xi\rangle$, $\xi \in \mathfrak g$ (a result
which will be proved in Lemma~\ref{KosLem}). The rest of this point is a 
consequence of Lemma~\ref{L31inv}. To prove the third point observe
that $C$ is a principal fibre bundle over the reduced space. Choose a
locally finite open covering $(U_\alpha)_{\alpha\in I}$ of 
$M_{\rm red}$ over which 
the bundle is trivial, choose local sections
$\sigma_\alpha:U_\alpha \to C$ of this bundle, and a smooth
partition of unity $(\psi_\alpha)_{\alpha\in I}$ subordinate to the
covering. Then the map
$\sigma(\phi):=\sum_{\alpha\in I}\psi_\alpha(\sigma_\alpha^*\phi)$ will
do the job. In case the fibres are compact one may also integrate out the
fibres with respect to some density to get the desired $\sigma$.
The last part is a direct computation.
\end{proof}

Note that for non-proper group actions it is in general no
longer true that each smooth complex-valued function on the
reduced phase space is induced by a globally $G$-invariant
function as the following example shows:

Let $M$ be the cotangent bundle of the two-torus $T^2$ minus
its zero-section. It is diffeomorphic to the Cartesian
product $T^2\times (\mathbb R^2\setminus \{0\})$. Let $G =
\mathbb R$ whence its Lie algebra and its dual are
canonically isomorphic to $\mathbb R$. Let $J: M\to \mathbb
R$ be the function $J(z_1,z_2,p_1,p_2) :=
\frac{1}{2}(p_1^2-p_2^2)$. 
It is easy to see that every value of $J$ is
regular, but that the orbit space of the corresponding constraint
surface is a smooth Hausdorff manifold if and only if the
value is equal to $0$: in this case the two angular
frequencies $\partial J/\partial p_1$ and $\partial
J/\partial p_2$ are always rationally dependent so that the
reduced phase space is symplectomorphic to two copies of the cotangent
bundle of the unit circle minus the zero-section. In all the other
cases a generic $G$-orbit is not closed, but has a closure
diffeomorphic to $T^2$. From this it easily follows that every globally
$G$-invariant smooth complex-valued function $f$ on $M$ is
of the general form $f(z_1,z_2,p_1,p_2) = \phi(p_1,p_2)$.
Hence its restriction to the constraint surface $C :=
J^{-1}(\{0\})$ does clearly only induce those functions on
$T^* S^1\setminus S^1$ which are invariant under
$U(1)$-rotations. In this example the difference between the
so-called `strongly invariant functions' (the elements of
$C^\infty(M)^G$) and the `weakly invariant functions' (the
elements of $\mathcal B_C$) becomes crucial.

Finally we should like to mention that the Marsden-Weinstein
reduction for Hamiltonian $G$-spaces for a
\emph{non-zero value} $\mu$ of the momentum map 
can be reduced to the above case by adding a suitable
coadjoint orbit, see e.g.~\cite[p.~194, Thm. 26.6]{GS84}.
Note also that the $G$-action on the extended system is proper if the
original action was proper.

\vspace{0.5cm}

We shall give another description of the space 
$C^\infty (M_{\rm red})$ as the zero group of BRST cohomology in the
next section. To this end let us recall briefly some results on Koszul 
and Chevalley-Eilenberg cohomology related to phase space reduction.

Let $V$ be an $n$-dimensional real vector space, $V^*$ its
dual, and $J: M \to V^*$ a smooth map such that $0$ is a
regular value of $J$ and the constraint surface 
$C := J^{-1}(\{0\})$ is nonempty. Then $C$ is a regular
sub-manifold of codimension $n$ of $M$. Define the 
\emph{ideal generated by $J$}, ${\mathcal I}(J)$, as the
ideal of $C^\infty(M)$ spanned by all
functions of the form $f\langle J, \xi\rangle$ where 
$\xi\in V, f\in C^\infty(M)$.
Note that this definition also makes
sense for any smooth map $M\rightarrow V^*$.
Denote by $\Grass V$ the Gra{\ss}mann algebra over $V$ and consider
the tensor product $\Grass V \otimes C^\infty(M)$. Let
\begin{equation}
\label{KoszulDef}
    \kos:
    \Grass^\bullet V \otimes C^\infty(M)
    \to
    \Grass^{\bullet-1} V \otimes C^\infty(M),
    \quad
    a \mapsto i(J)a
\end{equation}
denote the \emph{Koszul boundary operator} associated to $C$ and
$J$. Here $i(J)$ means the left insertion 
(the standard interior product) of $J$. We shall sometimes 
write $\kos_i$, $1\leq i \leq n$, for the restriction of $\kos$ to 
$\Grass^i V\otimes C^\infty(M)$. 
The pair $(\Grass V\otimes C^\infty(M),\kos)$ becomes a chain complex
(as $\kos^2=0$). For a
regular constraint surface this complex is known to be
acyclic, which can be seen using an augmentation (see
e.g. \cite[Def.~6.5, p.~339]{Jac89}): 
Let $(C^\infty(C)\oplus (\Grass V\otimes C^\infty(M)),\hat\kos)$ 
be the \emph{augmented Koszul complex}, 
where $\hat\kos$ is defined by 
$\hat{\kos}_i:=\kos_i$ for $1\leq i\leq n$ and 
$\hat{\kos}_{0} :=\iota^*$ (the \emph{augmentation}). 
\begin{lemma}
\label{KosLem}
With the above notations suppose that the constraint surface is a regular
value of the map $J$. Then there is a chain homotopy for the augmented 
complex: More precisely there is a linear map 
$\hat{h}$ with components 
$\hat h_{-1} = \prol: C^\infty (C) \to C^\infty (M)$, 
$\hat h_i = h_i: \Grass^i V \otimes C^\infty (M) \to 
\Grass^{i+1} V \otimes C^\infty (M)$, for $0 \le i \le n$, such that 
$\hat h \hat \kos + \hat \kos \hat h = \id$. Moreover, we can choose
$h_0$ such that
\begin{equation}
      \label{WichtigeGleichung}
          h_0 \prol = 0.
\end{equation}
In particular, the vanishing ideal ${\mathcal I}_C$ is equal to
the space of Koszul-0-boundaries which in turn is equal to the
ideal generated by $J$, ${\mathcal I}(J)$.
\end{lemma}
The additional technical equation (\ref{WichtigeGleichung}) will
become rather useful for the quantum deformation of all this in
Section \ref{BRSTCohoSec}. Again, for (proper) group actions there is
an equivariant analogue:
\begin{lemma}
\label{KosLemInv}
Under the above circumstances suppose in addition that a Lie
group $G$ acts on $M$ leaving invariant $C$. Suppose
furthermore that there is a representation of $G$ in $V$
such that $J$ is an equivariant map (with respect to the
contragredient representation of $G$ on $V^*$). Then all the maps
$\hat{\kos}_i$
are equivariant with respect to the natural action of $G$ on
$\Grass V \otimes C^\infty(M)$ and $C^\infty(C)$. Moreover, the chain 
homotopy $\hat{h}$ can in addition be chosen
to be an equivariant map in case the group action is proper.
\end{lemma}
\begin{proof}
Again we shall only be treating the $G$-invariant case and
use the notation of the proof of Lemma~\ref{L31inv}.
As in \cite[p.~9--10]{FK92} we shall first construct the chain homotopies
separately on a $G$-invariant open neighbourhood $U$ of $C$ and
on an open $G$-invariant set $W$ not meeting $C$ such that $U\cup W=M$.
The overlap region $U\cap W$ has to be treated with care to ensure
equation (\ref{WichtigeGleichung}). 

1. Start with a $G$-invariant tubular neighbourhood $U'$ as constructed
  in Lemma~\ref{L31inv} with $G$-equivariant projection $\tilde{\tau}$. Using the fact
  that $C$ is a regular constraint surface and techniques analogous to the ones used in 
  Lemma~\ref{L31inv} to establish injectivity of $\Phi$ we can arrange things
  in such a way that the map
  $U'\rightarrow C\times V^*:u\mapsto (\tilde{\tau}(u),J(u))$ becomes a 
  $G$-equivariant diffeomorphism onto an open $G$-invariant neighbourhood $Z$ of $C$ 
  in $C\times V^*$ when restricted to a suitable $G$-invariant open neighbourhood
  $U\subset U'$ of $C$ in $M$ which we shall often identify with its image $Z$ 
  in the sequel.
  Shrinking $U$ if necessary in a $G$ invariant manner allows us to assume that
  for each point $(c,\mu)\in Z$ the interval $\{(c,t\mu)|t\in[0,1]\}$ is also contained 
  in $Z$. Using a basis $e_1,\ldots,e_n$ of $V$ and linear co-ordinates
  $\alpha_1,\ldots,\alpha_n$ on $V^*$ we define the map
  $h_U: \Grass^{\bullet} V \otimes C^\infty(U) \to  
\Grass^{\bullet+1} V \otimes C^\infty(U)$ by   
\begin{equation} 
\label{HomotNearC} 
   h_U(\phi) (c,\mu) := \sum_{i=1}^n e_i\wedge
         \int_0^1 t^{k}\frac{\partial \phi}{\partial \alpha_i}
         (c,t\mu) dt,
\end{equation}
where $\phi\in\Grass^k V\otimes C^\infty(U)$.
It is a routine
check (similar to the proof of the Poincar\'e Lemma upon noting
that $J$ is equal to the projection of the second factor in $C\times V^*$)
that $h_U$ is a chain homotopy for the restriction of the Koszul boundary
operator to $\Grass V\otimes C^\infty(U)$ and that 
$\tau_U^*\iota^* + \kos h_U = \id_{C^\infty(U)}$ where
$\tau_U:=\tilde{\tau}|_U$.

2. Let $W$ be the complement of the closure of the set of all those points
in $U$ whose second co-ordinate in $Z$ is multiplied by $1/2$. $W$ is a 
$G$-invariant open subset of $M$ such that $M=U\cup W$. Let $1=\psi_U+\psi_W$
be a $G$-invariant smooth partition of unity subordinate to that covering.
We shall show in the next subsection that there is $G$-equivariant smooth
map $\xi:W\rightarrow V$ with support in $W$ such that 
i) $\langle J,\xi \rangle = \psi_W$ and 
ii) $\xi|_{\textrm{supp}(\psi_U)\cap W}=-h_U(\psi_U)$. Defining
for each $\phi\in \Grass V\otimes
C^\infty(W)$ the map $h_W(\phi):= \xi\wedge \phi$ it is not hard to check
that $h(\phi) := \psi_U h_U(\phi|_U) + h_W(\phi|_W)$
is the desired chain homotopy satisfying (\ref{WichtigeGleichung}) for the
geometric $G$-equivariant prolongation map constructed in Lemma~\ref{L31inv}
using $\psi_U$.

3. On $U\cap W$ define $\hat{\xi}:=-h_U(\psi_U)$ which is clearly
$G$-invariant and satisfies i) and ii) above. Let $W'$ be the
$G$-invariant open set $W\setminus \textrm{supp}(\psi_U)$. In order to define
$\xi$ on $W'$ with property i) above we proceed as follows: by the properness
of the $G$-action the isotropy subgroup $G_x$ of each point $x\in W'$ is compact.
Using a $G_x$-invariant scalar product on $V^*$ it is easy to construct
a smooth $G$-equivariant map $\eta^{[x]}$ of the $G$-orbit through $x$ into $V^*$
satisfying $\langle J,\eta^{[x]}\rangle =1 \quad \forall x \in W'$. 
Again by the properness of the action
each orbit is closed, and upon using a locally finite system of 
sufficiently small $G$-invariant tubular neighbourhoods 
around each orbit with subordinate $G$-invariant smooth
partition of unity and $G$-equivariant prolongation maps
(see again Lemma~\ref{L31inv}) we can glue together the prolongations of the 
maps $\eta^{[x]}$ to a $G$-equivariant smooth map $\xi':W'\rightarrow V^*$ satisfying
property i). The map $\xi$ is obtained by glueing $\hat{\xi}$ on $U\cap W$ and
$\xi'$ on $W'$ by means of a $G$-invariant smooth partition of unity subordinate
to the covering of $W$ by $U\cap W$ and $W'$.
\end{proof}

\noindent Again we shall call chain homotopies constructed in the above proof
\emph{geometric} chain homotopies.

Let $\rho:\mathfrak g \to \Hom(Q,Q)$ be a representation
of the finite-dimensional Lie algebra $\mathfrak g$ in some
complex vector space $Q$. Recall the definition of the 
\emph{Chevalley-Eilenberg differential}
$\ce: \Grass^\bullet \mathfrak g^* \otimes Q 
\to \Grass^{\bullet+1} \mathfrak g^* \otimes Q$: 
let
$\alpha\otimes q\in\Grass^k \mathfrak g^*\otimes Q$ and
$\xi_1,\ldots, \xi_{k+1}\in \mathfrak g$, then
\begin{eqnarray}
   \ce(\alpha\otimes q)(\xi_1,\ldots,\xi_{k+1})
                     & := &
                      \sum_{1\le i<j \le k+1}(-1)^{i+j+1}
                        \alpha([\xi_i,\xi_j], \xi_1, \ldots,
                        \stackrel{i}{\wedge}, \ldots,
                        \stackrel{j}{\wedge}, \ldots,
                        \xi_{k+1})\otimes q  \nonumber \\
                     &    & +
                      \sum_{i=1}^{k+1}(-1)^{i+1}
                        \alpha(\xi_1, \ldots,
                        \stackrel{i}{\wedge},
                        \ldots,\xi_{k+1})
                        \otimes \rho(\xi_i)(q). \label{CEDef}
\end{eqnarray}
It is well-known that $\ce^2=0$. We shall denote by
$H^\bullet_\ceind(\mathfrak g,Q)$ the Chevalley-Eilenberg cohomology
of $\mathfrak g$ with values in the $\mathfrak g$-module $Q$ for the
quotient of the kernel of $\ce$ by the image of $\ce$. This space
clearly inherits the $\mathbb Z$-grading of the Gra{\ss}mann algebra over
$\mathfrak g^*$.

For computations we shall frequently use a basis $e_1,\ldots,e_n$
($n:=\dim\mathfrak g$) of $\mathfrak g$
and its dual base $e^1,\ldots,e^n$ of $\mathfrak g^*$. Denoting by
${f^c}_{ab}:=\langle e^c, [e_a,e_b]\rangle$ the structure constants of
$\mathfrak g$ we get the following short formula for $\ce$:
\begin{equation}
\label{CEShortDef}
    \ce(\alpha\otimes q) = -\frac{1}{2}\sum_{a,b,c}{f^c}_{ab}\,
                           e^a\wedge e^b \wedge i(e_c) \alpha \otimes q
                           +\sum_{a}e^a\wedge \alpha\otimes \rho(e_a)(q)
\end{equation}
Recall that the zeroth Chevalley-Eilenberg cohomology group 
$H_\ceind^0(\mathfrak g,Q)$ is always equal to
the space of $\mathfrak g$-\emph{invariants} 
$Q^{\mathfrak g} := \{q \in Q \; | \; \rho(\xi)q=0\}$.

The significance of the Chevalley Eilenberg differential becomes clear
by the following obvious characterisation of the space of smooth
complex-valued functions on the reduced phase space:
\begin{proposition} 
\label{ZeroCohomClass}
Let $(M, \omega, G, J, C)$ be a Hamiltonian $G$-space with regular
constraint surface such that in addition the connected Lie group $G$
acts properly and freely on the constraint surface $C$. 
Then the space $C^\infty(M_{\rm red})$ is
canonically isomorphic to the zeroth Chevalley-Eilenberg cohomology of
$\mathfrak g$ with values in $\mathfrak g$-module $C^\infty(C)$ (with
representation $\ClassLieC$, see  (\ref{KlassLieM})) via the
pull-back $\pi^*$. For a general Hamiltonian Lie algebra action the
space $H^0_\ceind (\mathfrak g,C^\infty(C))$ will be called the space
of classical invariants on the constraint surface $C$. 
\end{proposition}

%
%

\section {The classical BRST cohomology with augmentation}
\label {ClassBRSTSec}

In this section we recall the construction of the classical BRST cohomology
following \cite{KS87}. Throughout the section 
$(M, \omega, \mathfrak g, J)$ is a Hamiltonian $\mathfrak g$-space. 
The classical BRST complex is a double complex whose total complex is
also a differential graded Poisson algebra, i.e.~the total differential is
a super-derivation of a super Poisson structure. It is constructed as
follows:
\begin {enumerate}
\item The space of chains is the $\mathbb Z \times\mathbb Z$-graded vector space
      $\mathcal A^{\bullet, \bullet} :=
      \Grass^{\bullet}\mathfrak g^* \otimes \Grass^{\bullet}\mathfrak g
      \otimes C^{\infty}(M)$,
      where the gradings are called by convention \emph{ghost} and
      \emph{antighost degree} (following \cite [p.~191]{HT92} where
      these gradings are called `pure ghost number' and 
      `antighost number'). Moreover, $\mathcal A$ carries a natural
      $\mathbb Z_2$-graded vector space structure
      $\mathcal A = \Grass^{\mathrm {even}} (\mathfrak g^* \oplus \mathfrak g)
                    \otimes C^{\infty} (M) \oplus
                    \Grass^{\mathrm {odd}}(\mathfrak g^* \oplus \mathfrak g)
                    \otimes C^{\infty} (M)$
      and a $\mathbb Z$-grading
      $\mathcal A^{(n)} = \bigoplus_{n=k-l} \mathcal A^{k,l}$,
      where $n \in \mathbb Z$ is by convention called 
      \emph{ghost number}, see \cite [p.~191] {HT92} (also called
      `total degree' in \cite[p.~57]{KS87}). Using the
      $\wedge$-product of forms (which of course is graded in the
      standard way, i.e. 
      $(\alpha \otimes \xi) \wedge (\beta \otimes \eta)
       = (-1)^{kl} (\alpha \wedge \beta) \otimes (\xi \wedge \eta)$,
       where $\alpha, \beta \in \Grass \mathfrak g^*$, 
      $\xi, \eta \in \Grass \mathfrak g$, and the degrees of $\beta$
      and $\xi$ are $k$, $l$, respectively) and the pointwise product
      of functions, $\mathcal A$ becomes an associative,
      super-commutative algebra, graded with respect to all the above
      mentioned degrees. We shall 
      sometimes use the physicist's terminology to call elements of 
      $\Grass \mathfrak g^*$ ghosts and elements of 
      $\Grass \mathfrak g$ antighosts. 

\item The vertical differential is taken to be the standard
      Chevalley-Eilenberg differential
      $\ce: \mathcal A^{\bullet,\bullet} \to \mathcal A^{\bullet +1,\bullet}$
      of the Lie algebra cohomology of $\mathfrak g$ with respect to the
      $\mathfrak g$-module $\Grass \mathfrak g \otimes C^{\infty}(M)$
      where the representation is given by
      $\mathfrak g \ni \xi \mapsto \ad(\xi) \otimes \id 
                          + \id \otimes \{\langle J,\xi\rangle,\cdot\}$,
      see (\ref{CEDef}). Its cohomology is denoted by 
      $H^{\bullet}_\ceind (\mathcal A)$.

\item The horizontal differential
      $\kos: \mathcal A^{\bullet, \bullet} 
      \to \mathcal A^{\bullet, \bullet -1}$
      is defined to be the standard extension of the previously
      defined Koszul differential (\ref {KoszulDef}) to the complex 
      $\mathcal A$, i.e.
      $\kos (\alpha \otimes x \otimes F) 
      = (-1)^k \alpha \otimes i(J) (x \otimes F)$
      for $\alpha \in \Grass^k \mathfrak g^*$, 
      $x \in \Grass \mathfrak g$, and $F \in C^\infty (M)$. Its
      homology will be denoted by $H_\bullet^\kosind (\mathcal A)$.

\item It is easy to see that $\kos$ anti-commutes with $\ce$ and hence
      we form out of this double complex its total complex as
      follows. The total differential 
      $\mathcal D: \mathcal A^{(\bullet)} \to \mathcal A^{(\bullet+1)}$
      is taken to be
      \begin{equation}
      \label{ClassBRSTOpDef}
          \mathcal D := \ce + 2 \kos
      \end{equation}
      and is called the \emph {(classical) BRST operator}
      (this inessentially differs from \cite{KS87} where the BRST
      operator was defined to be $\ce + 2(-)^k \kos,$ where $k$ is the
      antighost degree). 

\item The algebra $\mathcal A$ has a natural super Poisson structure
      induced by the natural pairing of $\mathfrak g$ and 
      $\mathfrak g^*$ and the Poisson bracket on $M$.
      In order to describe this bracket we firstly
      recall the definition of the left and right insertion maps on
      $\Grass(\mathfrak g^* \oplus \mathfrak g)$. Let
      $\alpha \in \Grass^k(\mathfrak g^*)$, 
      $\xi \in \Grass^l (\mathfrak g)$, $\beta \in \mathfrak g^*$, and 
      $X \in \mathfrak g$. Then 
      $i(\beta) (\alpha \otimes \xi) 
      := (-1)^k \alpha \otimes (i(\beta)\xi)$ and $i(X) (\alpha\otimes\xi) 
      := (i(X)\alpha) \otimes \xi$. Moreover the \emph{right insertion}
      $j(v)$ for $v \in \mathfrak g^* \oplus \mathfrak g$ is defined
      by $j(v) a = -(-1)^m i(v) a$ where 
      $a \in \Grass^m (\mathfrak g^* \oplus \mathfrak g)$.
      Then we define the following endomorphisms $P$ and $P^*$ of
      $\Grass(\mathfrak g^* \oplus \mathfrak g) 
      \otimes \Grass (\mathfrak g^* \oplus \mathfrak g)$ (here the
      tensor product is not graded)
      \begin{equation}
      \label{PdualPDef}
          P := \sum_{a=1}^n j(e_a) \otimes i (e^a)
          \qquad
          \textrm{ and }
          \qquad
          P^* := \sum_{a=1}^n j(e^a) \otimes i(e_a),
      \end{equation}
      which clearly do not depend on the choice of the basis. Then the 
      super Poisson bracket on 
      $\Grass(\mathfrak g^* \oplus \mathfrak g)$ is defined by
      \begin{equation}
      \label{NatPairPoissonDef}
          \{a, b\} 
          = 2 \mu\circ (P + P^*) (a \otimes b)
          = 2 \sum_{c=1}^n \left(j(e_c) a \wedge i(e^c) b + 
             j(e^c) a \wedge i(e_c) b\right),
      \end{equation}
      where $\mu$ denotes the $\wedge$-product of 
      $\Grass(\mathfrak g^* \oplus \mathfrak g)$. The factor $2$ is by 
      convention to get the correct Clifford algebra later. Now we
      tensor this with the Poisson bracket on $M$ 
      \begin{equation}
      \label{SuperPoisson}
           \{a \otimes F, b \otimes G\}
           := a\wedge b\otimes \{F,G\} + \{a,b\}\otimes FG,
      \end{equation}
      for $a, b \in \Grass(\mathfrak g^* \oplus \mathfrak g)$ and
      $F, G \in C^{\infty}(M)$,
      to get a super Poisson bracket on $\mathcal A^{(\bullet)}$.
      Note that this super Poisson bracket is still $\mathbb Z_2$-graded,
      but no longer graded with respect to ghost and antighost degree
      separately. It is, however, graded with respect to ghost number:
      this can also be seen by the Hamiltonian form of the ghost
      number derivation which we shall give further down.

\item The classical BRST operator turns out to be a Hamiltonian
      super-derivation of the above Poisson structure: regarding the
      Lie bracket $[\cdot, \cdot]$ of $\mathfrak g$ in a canonical
      manner as an element in 
      $\Grass^2 \mathfrak g^* \otimes \mathfrak g \subset \mathcal A$,
      we define
      \[
           \Omega := -\frac{1}{2} [\cdot, \cdot].
      \]
      More precisely, this means
      $\Omega = -\frac{1}{4} \sum_{a,b,c=1}^n 
      {f^c}_{ab}\, e^a \wedge e^b \wedge e_c$
      in terms of a chosen basis. On the other hand we have the momentum map
      $J \in \mathfrak g^* \otimes C^{\infty}(M) \subset \mathcal A$. So
      \begin{equation}\label{TetaClass}
          \Theta := \Omega + J
      \end{equation}
      is an odd element in $\mathcal A^{(1)}$ which we call the (classical)
      \emph{BRST charge}. One can easily verify that
      \begin{equation}\label{eq:BRSTCharge}
         \mathcal D = \{\Theta, \cdot\}.
      \end{equation}

\item Let $\gamma$ be the identity endomorphism of $\mathfrak g$ regarded as an
      element of $\mathcal A^{1,1}$. This takes the form
      \begin{equation}\label{Ghamma}
          \gamma = \frac{1}{2} \sum_{a=1}^n e^a \wedge e_a
      \end{equation}
      in terms of the aforementioned basis and dual basis. Note that
      the ghost number grading is induced by the \emph{ghost number derivation}
      \begin{equation} \label{GhostNumber}
          \mathsf{Gh} := \{\gamma,\cdot\}.
      \end{equation}
\end {enumerate}
\begin{definition}
The differential graded Poisson algebra 
$(\mathcal A^{(\bullet)}, \mathcal D, \{\cdot,\cdot\})$ is called the
classical BRST algebra. Its cohomology group 
$\ker {\mathcal D} \big/ \image {\mathcal D}$ will be called the
classical BRST cohomology and will be denoted by 
$H^{(\bullet)}_\brsind (\mathcal A)$.
\end{definition}
\begin{lemma}\label{ClassAlg}
The classical BRST cohomology is equipped with a natural
$\mathbb Z$-graded super Poisson
structure induced by the super Poisson structure of the classical BRST
algebra: let $a,b\in {\mathcal A}$ such that
${\mathcal D}a=0={\mathcal D}b$; then for the corresponding cohomology
classes $[a],[b]\in H^{(\bullet)}_\brsind (\mathcal A)$ we have
$[a]\wedge [b]:=[a\wedge b]$ and $\{[a],[b]\}:=[\{a,b\}]$.
\end{lemma}
\begin{proof}
Since the classical BRST operator is a super-derivation it immediately
follows that its kernel is a sub-algebra of the super Poisson algebra
$\mathcal A$. Moreover ${\mathcal D}^2=0$ entails that the image of
$\mathcal D$ is a super Poisson ideal in the kernel. The grading is
inherited by the induced ghost number derivation. This proves the
lemma.
\end{proof}

We shall now suppose that $(M, \omega, \mathfrak g, J, C)$ is a
Hamiltonian $\mathfrak g$-space with regular constraint surface $C$.
In order to get more information about the classical BRST cohomology
and its relation to the constraint surface $C$ we shall extend the
augmentation of the previous section to the above double complex:

For any $\xi\in\mathfrak g$ let $\xi_M$ and $\xi_C$ be the (infinitesimal)
generator of the action of $\mathfrak g$ or $G$ on $M$ and $C$,
respectively. The representation $\ClassLieC$ uniquely defines the 
Chevalley-Eilenberg complex 
$(\Grass \mathfrak g^*\otimes C^\infty(C), \ce^c)$ 
(as, of course, the representation 
$\ClassLieM$ defines the Chevalley-Eilenberg operator
$\ce$ on $\Grass \mathfrak g^*\otimes C^\infty(M)$). Now extend the
restriction map $\iota^*:C^\infty(M)\rightarrow C^\infty(C)$ to the
Chevalley-Eilenberg complex $\Grass \mathfrak g^*\otimes C^\infty(M)$ by
$\iota^*(\alpha\otimes f) := (-1)^k \alpha \otimes \iota^*f$ where 
$\alpha \in \Grass^k \mathfrak g^*$ and $f\in C^\infty(M)$. 
Thanks to the identity
\begin{equation}\label{CIntertw}
      \ClassLieC(\xi) \iota^* = \iota^* \ClassLieM(\xi)
\end{equation}
it is clear that
\begin{equation}\label{CCEAnticom}
        \ce^c \iota^* = - \iota^* \ce.
\end{equation}
For every prolongation map $\prol$, we have the following
equation, which will become important for deformation:
\begin{equation}
  \label{ClassConstrRep}
     \ClassLieC(\xi) = \iota^*\ClassLieM(\xi)\prol 
\end{equation}

Let $\hat{\mathcal A}$ denote the \emph{augmented classical BRST complex}
$(\Grass \mathfrak g^*\otimes C^\infty(C))\oplus \mathcal A$, and denote by
$\hat{\mathcal D}$ the \emph{augmented classical BRST operator}
\begin{equation}
\label{AugClasOp}
     \hat{\mathcal D} := \ce^c + 2\iota^* + \mathcal D,
\end{equation}
where all the maps are defined to be zero on the domains on which they
were previously not defined. Clearly
\begin{equation}
\label{AugClasEq}
       \hat{\mathcal D}^2 = 0.
\end{equation}
Moreover, we need to extend the chain homotopies $h$ and
$\prol$ of the classical (augmented) Koszul complex (compare 
Section~\ref{GeomClassSec}, Lemma~\ref{KosLem}) to $\mathcal A$ and 
$\hat{\mathcal A}$ which is done in the usual way by 
$h(\alpha\wedge f):= (-1)^k\alpha\wedge hf$ and
$\prol(\alpha\wedge \phi):=(-1)^k\alpha\wedge \prol\phi$,
where $\alpha\in\Grass^k\mathfrak g^*$, 
$f\in\Grass\mathfrak g \otimes C^\infty(M)$, $\phi\in C^\infty(C)$.
We keep the notation $\hat{h}$ for $\prol+h$ on 
$\hat{\mathcal A}$. Moreover, let $\hat{\ce}$ denote the augmented
Chevalley-Eilenberg operator $\ce^c+\ce$.

The cohomology of the classical BRST complex can be
computed in terms of the Chevalley-Eilenberg cohomology on the
constraint surface: 
\begin{proposition}
\label{cohomcrack}
Let $(M, \omega, \mathfrak g, J, C)$ be a Hamiltonian 
$\mathfrak g$-space with regular constraint surface.
With the above notations and definitions we have:
\begin{enumerate}
\item The following map
      \begin{equation}
      \label{ClassModHomDef}
           \hat{h}' := 
           \frac{1}{2}\hat{h} \left(\id +
           \frac{1}{2}(\hat{\ce}\hat{h}+\hat{h}\hat{\ce})\right)^{-1}
      \end{equation}
      is a chain homotopy for the augmented complex, i.e.
      $\hat{\mathcal D}\hat{h}' + \hat{h}'\hat{\mathcal D} = \id$.
      The linear map
      \begin{equation}
      \label{ClassKosIso}
          \Psi: H^{(\bullet)}_\brsind (\mathcal A)
          \to
          H^{\bullet}_\ceind(\mathfrak g,C^\infty(C)):  \quad
          [a] \mapsto [\iota^*a]
      \end{equation}
      is an isomorphism with the following inverse
      (where $[c]\in H^{\bullet}_\ceind(\mathfrak g,C^\infty(C))$):
      \begin{equation}
      \label{ClassKosInv}
          \Psi^{-1}: [c] \mapsto \left[2\hat{h}'c \right]
                         = \sum_{k=0}^{n}\left(-\frac{1}{2}\right)^k 
                         \left[(h\ce)^k\prol c\right].
      \end{equation}
\item The isomorphism $\Psi$ turns the Chevalley-Eilenberg cohomology
      on the constraint surface, 
      $H^{\bullet}_\ceind(\mathfrak g, C^\infty(C))$, into a 
      $\mathbb Z$-graded super Poisson algebra. The exterior
      multiplication and the super Poisson bracket take the following
      form for $c_1, c_2 \in \Grass \mathfrak g^*\otimes C^\infty(C)$:
      \begin{equation} 
      \label{ConstrCWedge}
          [c_1]\wedge [c_2] := [c_1\wedge c_2], 
      \end{equation}
      \begin{equation}
      \label{ConstCPoisson}
          \left\{[c_1],[c_2]\right\}   
          := \left[\iota^*\{\prol c_1,\prol c_2\}\right]
             -\frac{1}{2}
             \left[\iota^*\{\prol c_1, h_0 \ce \prol c_2\}\right]  
             -\frac{1}{2}
             \left[\iota^*\{h_0\ce\prol c_1,\prol c_2\}\right].
      \end{equation}
      In particular for the classical invariants, i.e.~the elements of 
      $H^0_\ceind(\mathfrak g, C^\infty(C))$, the Poisson bracket
      reduces to 
      \begin{equation} 
      \label{ClassPBra}
          \left\{[c_1],[c_2]\right\} = 
          \left[\iota^*\{\prol c_1,\prol c_2\}\right].
      \end{equation}
\end{enumerate}
\end{proposition}
\begin{proof}
For the first assertion, note that the chain homotopy
equation $\hat h \hat \kos + \hat \kos \hat h = \id$  
of the `pure' augmented Koszul complex (see Lemma \ref{KosLem}) still
holds on the augmented BRST complex. Hence
$\hat{\mathcal D}\hat{h}+\hat{h}\hat{\mathcal D}
= 2\id + \hat{\ce}\hat{h}+\hat{h}\hat{\ce}$.
Since $\hat{\ce}\hat{h}+\hat{h}\hat{\ce}$ is obviously nilpotent of
order at most $n+1$ the right hand side is invertible.
Thanks to $\hat{\mathcal D}^2 =0$ the map 
$\hat\delta \hat h + \hat h \hat \delta$ commutes with 
$\hat{\mathcal D}$. 
Thus $\hat h'$ is well-defined and a chain homotopy.
The fact that the map (\ref{ClassKosIso}) is well-defined
and bijective is shown using the chain homotopy property
and the equations $2\iota^*\hat{h}'\iota^*=\iota^*$ and
${\mathcal D}\hat{h}'\iota^*+\hat{h}'\ce^c \iota^*=0$ which
straightforwardly follow from the definitions by pure diagram chase. 
The simplification of
$\Psi^{-1}$ is a simple consequence of $c$ being closed with respect to
$\ce^c$. The second assertion is clear by Lemma \ref{ClassAlg}. The
form of the exterior multiplication is quickly computed using the fact
that $\mathcal A$ is $\mathbb Z \times \mathbb Z$-graded as a
super-commutative algebra. Note that the exterior multiplication of
the classes is well-defined since the representation $\ClassLieC$ is a
derivation on the $\mathfrak g$-module $C^\infty(C)$.
\end{proof}

We shall finally discuss the case when the Hamiltonian 
$\mathfrak g$-space comes from a Hamiltonian $G$-space where the
connected Lie group $G$ acts properly on $M$ and freely on $C$: We can
now choose the chain homotopies equivariant under the $G$-action 
(see Lemma~\ref{KosLemInv}).
Consequently 
\begin{equation}
  \label{ClassProperAnti}
   \hat{\ce}\hat{h}+\hat{h}\hat{\ce}=0,
\end{equation}
and the formula for the super Poisson bracket (\ref{ConstCPoisson})
simplifies to (\ref{ClassPBra}) for all elements of the
Chevalley-Eilenberg cohomology.
Note that the BRST cohomology is computing the so-called cohomology
along the leaves or $G$-orbits, see e.g. \cite[Thm.~3.8, p.~53]{Dub87} 
for a more precise statement. Furthermore a formula for the Poisson
bracket of two functions $\phi_1,\phi_2\in C^\infty(M_{\rm red})$ is
easily written down upon using a suitable left inverse $\sigma$ of the
pull-back with the projection $\pi:C\rightarrow M_{\rm red}$ 
(cf. Proposition~\ref{ReducClas}, \emph{iii.)}) and yields the same
result as in (\ref{PoissonUnten}).

%
%

\section {The quantum BRST operator in Deformation Quantization}
\label {QuantBRSTSec}

In this section the quantum BRST algebra will be defined and we shall
describe some operators associated with it. Throughout this section
$(M,*,\mathfrak g,\qmm)$ will be a fixed Hamiltonian quantum 
$\mathfrak g$-space. Although our aim is to describe the case of a
regular constraint surface (the irreducible case) our definitions
and results in this section do not need this restriction, neither do
we need the fact that $M$ is symplectic.

The underlying vector space for the quantum BRST algebra 
is the $\mathbb C[[\lambda]]$-module $\mathcal A[[\lambda]]$ of formal
power series with values in $\mathcal A$ endowed with its ghost and
antighost gradings as defined in the previous section. Moreover,
$\mathcal A[[\lambda]]$ inherits the ghost number grading from
$\mathcal A$ as well as the $\mathbb Z_2$-grading in even and odd
elements. Using the natural pairing of $\mathfrak g$ and 
$\mathfrak g^*$ as inner product we can define a one parameter family
of equivalent products $\kcliff$ for the Gra{\ss}mann part of
$\mathcal A[[\lambda]]$ indexed by a parameter $\kappa \in [0,1]$:
\begin{equation}
\label{KappaCliffProdDef}
    \alpha \kcliff \beta 
    = \mu \circ \eu^{2\im\lambda (\kappa P + (1-\kappa)P^*)}
      \alpha \otimes \beta,
\end{equation}
where the operators $P$ and $P^*$ are given as in (\ref{PdualPDef}).
This multiplication is known to be associative since left and right
insertions are anti-commuting super-derivations 
(see e.g.~\cite[Prop.~2.1]{Bor96} for a proof). 
Moreover, note that all the products $\kcliff$ satisfy the Clifford
relation on the one-forms for a multiple of the quadratic form on
$\mathfrak g^* \oplus \mathfrak g$ defined by the natural pairing.
All the multiplications $\kcliff$ are formal associative
deformations of the $\wedge$-product with first order super-commutator
being $\im\lambda$ times the super Poisson bracket 
(\ref{NatPairPoissonDef}). Taking now one fixed $\kappa$-ordered
product $\kcliff$ for the Gra{\ss}mann part and the  
quantum covariant star product $*$ for the functions we obtain an
associative product for $\mathcal A[[\lambda]]$ by tensoring these
algebra structures which we shall denoted by $\kstar$: more precisely,
for $\alpha \otimes F, \beta \otimes G \in \mathcal A[[\lambda]]$ with 
$\alpha,\beta \in 
\Grass (\mathfrak g^* \oplus \mathfrak g)[[\lambda]]$ 
and $F,G \in C^\infty (M)[[\lambda]]$ we set
\begin{equation}
\label{KappaOrderedProduct}
    (\alpha \otimes F) \kstar (\beta \otimes G )
    = (\alpha \kcliff \beta) \otimes (F * G).
\end{equation}
Then all $\kstar$ are equivalent by means of the equivalence
transformation $S_\kappa = \exp(2\im\kappa\lambda\Delta)$
where $\Delta$ is defined by
\begin{equation}
\label{SuperLaplacian}
    \Delta := \sum_{a=1}^n i(e_a)i(e^a) = - \sum_{a=1}^n j(e_a)j(e^a).
\end{equation}
Clearly $\Delta$ does not depend on the choice of the basis
and the equivalence 
\begin{equation}
\label{KappaEquivStd}
    S_\kappa (a \kcliff b) = S_\kappa a \std S_\kappa b
\end{equation}
is a straightforward computation where 
$a,b \in \Grass (\mathfrak g^* \oplus \mathfrak g)$.
It can easily be seen that $\kstar$ is still
$\mathbb Z_2$-graded but $\kstar$ is no longer graded with respect to
the ghost and antighost degree separately. It will, however, be graded
with respect to ghost number which can be seen by means of the
classical ghost number derivation which turns out to be a derivation
of $\kstar$. Taking super-commutators with respect to $\kstar$ will be
denoted by $\adk$. In the particular cases of the Weyl 
(i.e.~$\kappa = \frac{1}{2}$) and standard ordered 
(i.e.~$\kappa = 0$) multiplications
we shall use the notation $\wpr$ and $\spr$ for 
the star products, respectively, and $\adw$ and $\ads$ for the 
super-commutators, respectively.

Now we should like to turn $\mathcal A[[\lambda]]$ into a 
\emph{graded differential algebra} by defining a cochain complex
whose differential is an odd $\mathbb C[[\lambda]]$-linear 
left super-derivation raising the ghost number by one. This is
achieved in the following way: We define the 
\emph{quantum Weyl ordered BRST charge} to be  
\begin{equation}
\label{WeylBRSTChargeDef}
    \wbrsc := \Omega + \qmm
\end{equation}
which is odd and of ghost number $1$. Associated to $\wbrsc$ is the 
\emph{Weyl ordered BRST operator}:
\begin{equation}
\label{WeylBRSTOp}
    \wbrso := \frac{1}{\im \lambda} \adw (\wbrsc).
\end{equation}
Note that $\wbrso$ is well defined because $\wpr$ is super-commutative 
in  $0$-th order of $\lambda$. Moreover $\wbrso$ is odd. Under the
equivalence transformations $S_\kappa$, the BRST charge and the BRST
operator transform according to:  
\begin{equation}
\label{KappaBRST}
    \kbrso := S^{-1}_{\kappa-\frac{1}{2}}\circ\wbrso \circ
    S_{\kappa-\frac{1}{2}},
\end{equation}
\begin{equation}
\label{KappaBRSTCharge}
    \kbrsc := S^{-1}_{\kappa-\frac{1}{2}}(\wbrsc),
\end{equation}
such that $\kbrso = (\im\lambda)^{-1}\adk (\kbrsc)$. The standard
ordered quantum BRST operator (where $\sbrsc = \Theta_0$)
\begin{equation}
\label{StandBRSTOp}
    \sbrso := \frac{1}{\im \lambda} \ads (\sbrsc)
\end{equation}
will turn out to be of major importance in the sequel.
\begin{definition} 
\label{BRST-algebra}
The triple $(\mathcal A^{(\bullet)}[[\lambda]],\kstar, \kbrso)$ is
defined  to be the $\kappa$-ordered quantum BRST algebra.
\end{definition}
\begin{lemma} \label{BRSTDifferentialLem}
For all $\kappa \in [0,1]$ we have:
\begin{enumerate} 
\item $\kbrsc = \Omega + \qmm + \im\lambda(1-2\kappa) \chi$
      where 
      $\chi \in \mathfrak g^* \subset \mathcal A^{1,0}[[\lambda]]$
      defined by $\chi (\xi) = \frac{1}{2} \tr (\ad (\xi))$ for 
      $\xi \in \mathfrak g$ is the trace form of $\mathfrak g$.
\item $\kbrsc \kstar \kbrsc = 0$.  
\item The classical ghost number derivation 
      $\mathsf{Gh} = \{\gamma, \cdot\}$, see 
      (\ref{Ghamma}, \ref{GhostNumber}), is equal to 
      $\frac{1}{\im\lambda} \adk (\gamma)$
      for all $\kappa\in [0,1]$ and induces the ghost number grading: 
      $\mathsf{Gh} (\Phi) = n\Phi \iff \Phi \in \mathcal A^{(n)}[[\lambda]]$
      for all $\Phi \in \mathcal A[[\lambda]]$.
\item $[\mathsf{Gh},\kbrso]=\kbrso$, hence $\kbrso$ raises the ghost number
      by one.
\end{enumerate}
\end{lemma}
\begin{proof}
The first part is a simple computation using the easily verified fact
that $\Delta \Omega = \chi$. For the second part, thanks to 
(\ref{KappaBRSTCharge}) it suffices to show 
$\wbrsc \wpr \wbrsc = 
\Omega \weyl \Omega + \Omega \wpr \qmm + \qmm \wpr \Omega 
+ \qmm\wpr\qmm = 0$. 
We compute all these terms: Firstly notice that 
$P^2 \Omega\otimes\Omega = 0 = (P^*)^2 \Omega \otimes \Omega$.
Since $PP^* = P^*P$ we only have to compute the orders $\lambda^0$, 
$\lambda^1$, and $\lambda^2$ of $\Omega \weyl \Omega$. The zeroth
order is trivially zero since $\Omega$ is odd, the first order
vanishes by use of the Jacobi identity for the Lie bracket of
$\mathfrak g$. Finally the second order is shown to  
vanish by direct computation of $\mu \circ (P^*P \Omega \otimes \Omega) = 0$. 
Thus $\Omega \weyl \Omega = 0$. Writing $\qmm = \sum_{c} e^c \otimes \qmm_c$
with $\qmm_c \in C^\infty (M)[[\lambda]]$ using a basis of $\mathfrak g$ 
we compute the anti-commutator
\[
    \Omega \wpr \qmm + \qmm \wpr \Omega 
    = \sum_c (\Omega \weyl e^c + e^c \weyl \Omega) \otimes \qmm_c
    = \sum_c \im\lambda (j(e^c) \Omega + i(e^c) \Omega) \otimes \qmm_c 
    = -\frac{\im\lambda}{2} \sum_{a, b, c} {f^c}_{ab}\,
      e^a \wedge e^b \otimes \qmm_c.
\]
Finally we have due to the quantum covariance of the star product $*$
\[
    \qmm \wpr \qmm = \sum_{a, b} e^a \wedge e^b \otimes \qmm_a * \qmm_b
                   = \frac{1}{2} \sum_{a, b} e^a \wedge e^b 
                     \otimes (\qmm_a * \qmm_b - \qmm_b *\qmm_a)
                   = \frac{\im\lambda}{2} \sum_{a, b, c} 
                     {f^c}_{ab}\, e^a \wedge e^b \otimes \qmm_c.
\]
Thus we have proved $\kbrsc \kpr \kbrsc = 0$.
The third part is again easily verified. Finally the fourth part 
follows from the third part: 
$[\mathsf{Gh},\kbrso] = [\mathsf{Gh}, \frac{1}{\im\lambda} \adk (\kbrsc)]
= \frac{1}{\im\lambda}\ad_{\kpr}(\mathsf{Gh}\kbrsc) = \kbrso$ since clearly
$\mathsf {Gh} \kbrsc = \kbrsc$.
\end{proof}
\begin{proposition} \label{dgA}
The $\kappa$-ordered BRST algebra is a differential graded algebra
over $\mathbb C [[\lambda]]$ for all $\kappa \in [0,1]$. All
of these are isomorphic as differential graded algebras via $S_\kappa$ and are 
formal deformations of the classical BRST algebra (their `classical limit').
\end{proposition}
\begin{proof}
It remains to compute the classical limit which is straightforward.
\end{proof}

In order to prepare Theorem \ref{exBRST}, in which the explicit form of
the  quantum BRST operator will be given, we make the following
definitions:
\begin{definition}
\label{defqrcu}
We define the following $\mathbb C[[\lambda]]$-linear endomorphisms of 
$\mathcal A^{\bullet,\bullet}[[\lambda]]$: For all 
$1\leq k,l\leq n$, $\alpha \in \Grass^k \mathfrak g^*$, 
$\xi = \xi_1 \wedge \cdots \wedge \xi_l \in \Grass^l \mathfrak g$, and
$F \in C^\infty (M)[[\lambda]]$ we define
\begin{enumerate}
\item $\quop:  \mathcal A^{\bullet,\bullet}[[\lambda]] 
       \to \mathcal A^{\bullet,\bullet -1}[[\lambda]]$ by
      \begin{equation}
      \label{QuOpDef}
          \quop (\alpha \wedge \xi \otimes F) :=
          (-1)^k \alpha \wedge
          \sum_{i<j} (-1)^{i+j-1} [\xi_i, \xi_j] \wedge \xi_1 
          \wedge \cdots 
          \stackrel{i}{\wedge} \cdots \stackrel{j}{\wedge} \cdots
          \wedge \xi_l \otimes F.
      \end{equation}

\item $\smu, \amu: \mathcal A^{\bullet,\bullet}[[\lambda]] 
      \to \mathcal A^{\bullet,\bullet -1}[[\lambda]]$ by
      \begin{equation}
      \label{MsMaDef}
      \begin{array} {rcl}
          \smu (\alpha \wedge \xi \otimes F) 
          & := & (-1)^k \sum_{a} 
                 \alpha \wedge i(e^a) \xi \otimes F * \qmm_a \\
          \amu (\alpha \wedge \xi \otimes F)
          & := & (-1)^k \sum_{a}
                 \alpha \wedge i(e^a) \xi \otimes \qmm_a * F,
      \end{array}
      \end{equation}
      where $\qmm = \sum_{a} e^a \otimes \qmm_a$ with 
      $\qmm_a \in C^\infty (M)[[\lambda]]$.

\item $\cuop: \mathcal A^{\bullet,\bullet}[[\lambda]]
      \to \mathcal A^{\bullet -1,\bullet -2}[[\lambda]]$ by
      \begin{equation}
      \label{COpDef}
          \cuop (\alpha \wedge \xi \otimes F) :=
          \sum_{i<j} (-1)^{i+j-1} i ([\xi_i, \xi_j])\alpha 
          \wedge \xi_1 \wedge \cdots 
          \stackrel{i}{\wedge} \cdots \stackrel{j}{\wedge} \cdots
          \wedge \xi_l \otimes F.           
      \end{equation}

\item $\uni: \mathcal A^{\bullet,\bullet}[[\lambda]]
      \to \mathcal A^{\bullet,\bullet -1}[[\lambda]]$ by
      \begin{equation}
      \label{UniDef}
          \uni := 
          \frac{1}{\im\lambda} \ad_\kappa (\chi) = \{ \chi, \cdot\}, 
      \end{equation}     
      where $\chi \in \mathfrak g^*$ is the trace form as in 
      Lemma~\ref{BRSTDifferentialLem}. 
\end{enumerate}
\end{definition}
Clearly it is sufficient to specify these operators only on the above
factorising elements. Note that $\uni$ does not depend on $\kappa$ but
consists only of the lowest order term. For later use we give for
these operators the following expressions in terms of a basis of
$\mathfrak g$ and $\mathfrak g^*$.
\begin{lemma}
\label{QCUBasisLem}
In terms of a basis $e_1, \ldots, e_n$ of $\mathfrak g$ and the dual
basis $e^1, \ldots, e^n$ of $\mathfrak g^*$ we have the following
expressions:
\begin{eqnarray}
    \quop 
    & = & - \frac{1}{2} \sum_{a,b,c} 
          {f^c}_{ab}\, e_c \wedge i(e^a)i(e^b) \label{QuOpBasis} \\
    \cuop
    & = & - \frac{1}{2} \sum_{a,b,c}
          {f^c}_{ab}\, i(e_c)i(e^a)i(e^b) \label{CuOpBasis} \\
    \uni
    & = & \sum_{a,b} {f^b}_{ab}\, i(e^a) \label{UniBasis}
\end{eqnarray}
\end{lemma}
\begin{proof}
This is a straightforward computation. 
\end{proof}

Using the $\mathfrak g$-representation $\QuantLieM$ as in
(\ref{QuantLieMDef}) on $C^\infty (M)[[\lambda]]$ the map   
\begin{equation}
\label{BRSTQuantRep}
    \mathfrak g \ni \xi \mapsto 
    \ad (\xi) \otimes \id + \id \otimes \QuantLieM (\xi)
\end{equation}
turns $\Grass \mathfrak g \otimes C^\infty (M)[[\lambda]]$ into a
$\mathfrak g$-module. This motivates the following definition:
\begin{definition}\label{defqds}
\begin{enumerate}
\item The Chevalley-Eilenberg differential 
      $\mathcal A^{\bullet,\bullet} [[\lambda]] \to 
      \mathcal A^{\bullet+1,\bullet} [[\lambda]]$ with respect to the 
      representation (\ref{BRSTQuantRep}) is denoted by $\qce$ and
      called the quantised Chevalley-Eilenberg differential.
\item The operator $\qkos: \mathcal A^{\bullet, \bullet}[[\lambda]]
      \to \mathcal A^{\bullet, \bullet-1}[[\lambda]]$ defined by
      \begin{equation}
      \label{quantKosDef}
          \qkos := \smu + \im\lambda(\frac{1}{2} \uni -\quop)
      \end{equation}
      is called the quantised Koszul differential. We shall frequently
      write $\qkos_i$ for its restriction to $\mathcal A^{\bullet,i} [[\lambda]]$,
      $(0\leq i\leq n)$.
\end{enumerate}
\end{definition}
Note that the classical limits of the operators $\qce$ resp. $\qkos$
are indeed $\ce$ resp. $\kos$. In the following theorem we shall
relate the standard ordered BRST operator $\sbrso$ to the more
concrete operators $\qce$ and $\qkos$.
\begin{theorem} 
\label{StandardBRSTDoubleTheo}
The standard ordered BRST operator satisfies
\begin{equation}
\label{StandardBRSTOperator}
    \sbrso = \qce + 2 \qkos
\end{equation} 
and thus defines a double complex, i.e.
\begin{equation}
\label{QuantumComplex}
    \qce^2 = 0, \quad \qkos^2 = 0, \quad \qce \qkos + \qkos\qce = 0.
\end{equation}
\end{theorem}
\begin{proof}
Clearly (\ref{QuantumComplex}) follows from
(\ref{StandardBRSTOperator}) and $\sbrso^2=0$ since $\qce$
resp. $\qkos$ are homogeneous of bidegree $(1,0)$ resp. $(0,-1)$. Thus 
we have to show (\ref{StandardBRSTOperator}) which is a
straightforward computation using 
$\sbrsc = \Omega + \qmm + \im\lambda\chi$. In terms of the basis we
obtain by a simple computation
\[
    \ads (\Omega) 
    = \im\lambda\left( \sum_{a,b,c} {f^c}_{ab}\,
      e^a \wedge e_c\wedge i(e^b)
          - \frac{1}{2} \sum_{a,b,c} {f^c}_{ab}\, 
          e^a \wedge e^b \wedge i (e_c) \right) 
    - \lambda^2 \sum_{a,b,c} 
    {f^c}_{ab}\, e_c \wedge i (e^a) i(e^b) 
\]
and similar $\ads (\qmm) = \sum_a e^a \wedge \ad_* (\qmm_a) + 2\im\lambda \smu$.
Together with (\ref{UniDef}), Lemma~\ref{QCUBasisLem}, and
(\ref{CEShortDef}) applied for $\qce$ the theorem follows.
\end{proof}

As we have seen $\sbrso$ splits into two super-commuting
differentials $\qce$ and $\qkos$, a fact which is only visible in this
clarity in the standard ordered case. In the general $\kappa$-ordered
case and in particular in the Weyl case a direct computation of
$\kbrso$ would be less useful. Hence we transform $\qce$ and $\qkos$
back via the equivalence transformation $S_\kappa$ to find the
differentials also in the $\kappa$-ordered case. We define
\begin{equation}
\label{kappaQCEQKosDef}
    \qce_\kappa := S^{-1}_\kappa \qce S_\kappa
    \quad
    \textrm { and }
    \quad
    \qkos_\kappa := S^{-1}_\kappa \qkos S_\kappa.
\end{equation}
It turns out that neither $\qce_\kappa$ nor
$\qkos_\kappa$ respect the ghost/antighost degree. To find a more
explicit form for $\qce_\kappa$ and $\qkos_\kappa$ we need the
following lemma:
\begin{lemma}
\label{kappaQCMUOperatorLem}
Let $\kappa \in [0,1]$ then 
$\Delta \quop - \quop \Delta = \cuop$, and
$\Delta$ commutes with $\cuop$, $\smu$, $\amu$, and $\uni$.
Hence $S^{-1}_\kappa \quop S_\kappa = \quop - 2\im\kappa\lambda\cuop$, 
and $S^{-1}_\kappa \smu S_\kappa = \smu$, 
$S^{-1}_\kappa \amu S_\kappa = \amu$, and 
$S^{-1}_\kappa \uni S_\kappa = \uni$. Finally
\begin{equation}
\label{kappaQCE}
    \begin{array}{c}
    \Delta \qce - \qce \Delta 
    = -2\quop - \frac{1}{\im\lambda} (\amu - \smu) + \uni, \\
    S^{-1}_\kappa \qce S_\kappa 
    = \qce + 4\im\kappa\lambda\quop + 2\kappa (\amu - \smu) 
    -2\im\kappa \lambda\uni + 4\kappa^2\lambda^2\cuop.
    \end{array}
\end{equation}
\end{lemma}
\begin{proof}    
Using the expressions in terms of a basis for these various operators
the above commutation relations follow from a tedious but straightforward
computation. 
\end{proof}

Using this lemma we can now state the following theorem which gives an 
explicit form for all $\kappa$-ordered BRST operators.
\begin{theorem}
\label{exBRST}
Let $\kappa \in [0,1]$ then the $\kappa$-ordered BRST operator
$\kbrso$ splits into two super-commuting
differentials $\kbrso = \qce_\kappa + 2 \qkos_\kappa$, where
\begin{eqnarray}
    \qce_\kappa 
    & = & \qce + 4 \im\kappa\lambda\quop - 2\kappa(\smu - \amu)
          - 2\im\kappa\lambda\uni + 4 \kappa^2\lambda^2\cuop, 
          \label{kappaQCEex} \\
    \qkos_\kappa
    & = & \qkos - 2 \kappa \lambda^2 \cuop = \smu +
          \frac{\im\lambda}{2} \uni - \im\lambda\quop -
          2\kappa\lambda^2\cuop, 
          \label{kappaQKosex} \\
    \kbrso    
    & = & \qce + 2\left((1-\kappa)\smu + \kappa \amu\right) 
          + 2\im\lambda(2\kappa-1)\quop - \im\lambda(2\kappa-1)\uni 
          - 4\kappa(1-\kappa) \lambda^2 \cuop.
          \label{kappaBRSTex}
\end{eqnarray}
\end{theorem}
\begin{proof}
This follows from the very definitions and the last lemma.
\end{proof}
\begin{remark}
In the standard ordered case the splitting of $\sbrso$ is rather
simple to find since in this case the two differentials are
homogeneous of bidegree $(1,0)$ resp. $(0,-1)$. Nevertheless for
physical reasons one is mostly interested in the Weyl ordered case
since only here the pointwise complex conjugation is a super
involution (a fact which we shall not need in this paper). 
But in the Weyl case the splitting is 
less obvious (and not compatible with the degrees) whence the
standard ordered case is a useful tool. In particular we have
\begin{equation}
\label{WeylCaseBRSTExplicit}
    \begin{array}{rcl}
    \wqce 
    & = & \qce + 2\im\lambda\quop - (\smu - \amu) 
          - \im\lambda\uni + \lambda^2 \cuop, \\
    \wqkos 
    & = & \smu + \frac{\im\lambda}{2}\uni - \im\lambda\quop 
          - \lambda^2\cuop, \\
    \wbrso
    & = & \qce + \smu + \amu - \lambda^2 \cuop.
    \end{array}
\end{equation}
\end{remark}

We shall denote the cohomology
$\ker \sbrso \big/ \image \sbrso$ of the quantised standard ordered
BRST differential by 
$\qcohom^{(\bullet)}_\brsind(\mathcal A [[\lambda]])$
and shall speak of the \emph{quantum BRST cohomology}. The cohomology
of $\qce$ will be denoted by 
$\qcohom_\ceind^{\bullet} (\mathcal A[[\lambda]])$, and the homology
of $\qkos$ by 
$\qcohom_\bullet^\kosind(\mathcal A[[\lambda]])$. 
Exactly as in the classical case we get the following quantum analogue of 
Lemma~\ref{ClassAlg}: 
\begin{lemma}\label{QuantAlg}
The quantum BRST cohomology is equipped with a natural
$\mathbb Z$-graded associative $\mathbb C[[\lambda]]$-bilinear
multiplication $\spr$ induced by the associative multiplication $\spr$
of the quantum BRST algebra ${\mathcal A}[[\lambda]]$: let 
$a, b \in \mathcal A$ such that
$\sbrso a=0=\sbrso b$; then for the corresponding cohomology
classes $[a],[b]\in H^{(\bullet)}_\brsind (\mathcal A [[\lambda]])$ we
have $[a]\spr [b]:=[a\spr b]$.

Moreover, the equivalence transformation $S_\kappa$ renders all the 
cohomologies of the $\kappa$-ordered BRST differentials canonically
isomorphic as associative algebras for all $\kappa \in [0,1]$.
\end{lemma}
The proof is completely analogous to the proof of Lemma \ref{ClassAlg}.

We shall now introduce a quantum analogue of the classical ideal
generated by $J$, 
${\mathcal I}(J)$
(see Section \ref{GeomClassSec} for a definition) and its normalizer
thereby generalising a notion also used by Gl\"o{\ss}ner in
\cite{Gloe98a}:  
\begin{definition} \label{QuantVanIdeal}
\begin{enumerate}
\item Let $\qideal (\qmm) := \qkos({\mathcal A}^{0,1}[[\lambda]])
      \subset{\mathcal A}^{0,0}[[\lambda]]$ be the quantum ideal of $\qmm$.
\item Let 
      $\qbideal (\qmm) := \left\{f \in C^\infty(M)[[\lambda]]
      \; \big| \; f * g - g * f \in \qideal (\qmm) 
      \quad \forall g \in \qideal (\qmm) \right\}$
      denote the quantum idealiser of $\qideal (\qmm)$.
\end{enumerate}
\end{definition}
\begin{proposition}\label{QIdealRep}
\begin{enumerate}
\item $\qideal(\qmm)$ is a left ideal of the algebra $(C^\infty(M)[[\lambda]],*)$.
      Moreover $\qideal(\qmm)$ is stable under the representation $\QuantLieM$,
      see (\ref{QuantLieMDef}).
\item $\qbideal(\qmm)$ is a sub-algebra of the algebra $(C^\infty(M)[[\lambda]],*)$
      containing $\qideal (\qmm)$ as a two-sided ideal.
\item The quotient $\qbideal (\qmm)/\qideal (\qmm)$ becomes an associative algebra
      in a canonical way.
 \end{enumerate}
\end{proposition}
\begin{proof}
First note that the operator $\quop$ vanishes on 
${\mathcal A}^{0,1}[[\lambda]]$, hence every element of 
$\qideal(\qmm)$ is equal to a sum of elements of the form 
$a*(\langle \qmm + \lambda \chi, \xi\rangle)$ where 
$a\in C^\infty(M)[[\lambda]]$ and $\xi \in \mathfrak g$ which
immediately shows that $\qideal (\qmm)$ is indeed a left ideal of the
algebra $C^\infty(M)[[\lambda]]$. Furthermore, $\qideal (\qmm)$ is 
stable under the above representation thanks to quantum covariance and
the fact that $\chi$ vanishes on commutators.
The second part is true by abstract algebra and follows from
associativity and the third part is clear.
\end{proof}

Note that $\qbideal (\qmm)$ is the idealiser of the left ideal
$\qideal(\qmm)$ in the sense of \cite[Eq.~(21), p.~199]{Jac89}.
In physics the ideal $\qideal (\qmm)$ may roughly be interpreted as
the space of all those `operators in the big unphysical Hilbert space'
vanishing on the smaller `physical Hilbert space', whereas
$\qbideal(\qmm)$ will then be the space of all those `operators in the
big unphysical Hilbert space' leaving invariant the `physical
subspace' whence the quotient algebra serves as the observable algebra
on the physical Hilbert space.

At the end of this section we should like to mention an important
algebra homomorphism relating the ghost-number zero part of the
quantum BRST cohomology to the above quotient algebra in Proposition
\ref{QIdealRep}: 
\begin{proposition} \label{BRSTDirac}
Let $(M,*,\mathfrak g,\qmm)$ be a Hamiltonian quantum 
$\mathfrak g$-space. Then there is a canonical homomorphism $L$ of
associative algebras 
\begin{equation}
\label{BRSTDiracHomo}
    L:\qcohom_\brsind^0(\mathcal A[[\lambda]])\rightarrow 
                            \qbideal (\qmm)/\qideal (\qmm),
\end{equation}
mapping each point $[\phi=\sum_{i=0}^n\phi_i] \in
\qcohom_\brsind^0(\mathcal A[[\lambda]])$ where $\phi_i\in
\mathcal A^{i,i}$ and $\sbrso\phi=0$ to $\phi_0 \bmod \qideal(\qmm)$.
\end{proposition}
\begin{proof}
Note that the equation $\sbrso\phi=0$ means in lowest order that
for $\phi_0$ there is $\phi_1$ such that $\qce \phi_0 +2\qkos_1\phi_1=0$.
Evaluating on $\mathfrak g$ this means that 
$\langle \qmm, \xi\rangle * \phi_0 - \phi_0 * \langle \qmm, \xi\rangle$
is contained in $\qideal (\qmm)$. Since $\chi(\xi)$ obviously commutes
with every $\phi_0$ and using associativity we see that $\phi_0$ has
to be in $\qbideal(\qmm)$. 
Moreover in case $\phi=\sbrso \phi'$ this means that
$\phi_0=\qkos_1\phi'_0$ where 
$\phi'_0\in \mathfrak g\otimes C^\infty(M)[[\lambda]]$ whence the BRST
quantum coboundaries are mapped to the quantum ideal of $\qmm$ which
shows that the map $L$ is well-defined. In order to see that $L$ is a
homomorphism of associative algebras one only has to observe that the
component in $\mathcal A^{0,0}$ of the multiplication $\phi\spr\psi$
(where again $\psi\in \mathcal A^{(0)}[[\lambda]]$ 
with $\psi=\sum_{i=0}^n\psi_i$, $\psi_i\in
\mathcal A^{i,i}$ and $\sbrso\psi=0$) is simply equal to
$\phi_0*\psi_0$ thanks to the fact that 
\[
    \mathcal A^{i,j}[[\lambda]] \spr \mathcal A^{k,l}[[\lambda]] 
    \subset 
    \sum_{r=0}^{\min\{j,k\}}\mathcal A^{i+k-r,j+l-r}[[\lambda]]
\]
which is particular for the standard ordered product $\spr$ on 
$\mathcal A[[\lambda]]$.
\end{proof}

%
%

\section {Computation of the quantum BRST Cohomology}
\label {BRSTCohoSec}

In Theorem \ref{exBRST} we have shown that the standard
ordered quantum BRST operator forms a double complex. 
As already mentioned both quantised differentials
$\qce$, $\qkos$ are deformations of the
corresponding classical operators $\ce$, $\kos$. This enables us to 
simplify the quantum BRST cohomology in exactly the same way as it was 
done for the classical cohomology in Section \ref{ClassBRSTSec}.

In order to get the full quantum analogue of Lemma
\ref{KosLem} and of Proposition \ref{cohomcrack} we still have to
`quantise' the augmentation, i.e.~to define a reasonable deformation
of the restriction map $\iota^*$. 
In the following proposition all the maps will be defined on the whole
quantum BRST algebra. Moreover, we shall suppose from now on that 
$(M, \omega)$ is symplectic and $(M, *, \mathfrak g, \qmm, C)$ is a
Hamiltonian quantum $\mathfrak g$-space with regular constraint
surface, and we shall write  
$\qideal_C$ for the quantum ideal generated by $\qmm$ now speaking of
it as the \emph{quantum vanishing ideal of $C$}. Likewise, its quantum
idealiser $\qbideal (\qmm)$ will now be denoted by $\qbideal_C$ to get
an analogy to the classical case (see Lemma~\ref{ReducClas}).
\begin{proposition} \label{QAug} 
Using the notation of Lemma \ref{KosLem} and of the two preceding
sections we define $B_1:= \frac{1}{\lambda}(\kos_1-\qkos_1)$ and,
using the chain homotopy $h_0$ (see Lemma \ref{KosLem}), the deformed
restriction map 
\begin{equation}\label{QuantRestrict}
    \qrestr
    := \iota^*(\id_{0}-\lambda B_1h_0)^{-1}
    = \sum_{r=0}^\infty \lambda^r \qrestr_r,
\end{equation}
where we have written
$\id_{0}:=\id_{\Grass \mathfrak g^*\otimes C^\infty(M)[[\lambda]]}$.
Then $\qrestr$ is the unique 
$\mathbb C [[\lambda]]$-linear map
$\Grass \mathfrak g^*\otimes C^\infty(M)[[\lambda]] \to
\Grass \mathfrak g^*\otimes C^\infty(C)[[\lambda]]$ 
which satisfies
\begin{eqnarray}\label{QuantAug}
    \qrestr_0        & = & \iota^*, \\
    \qrestr \qkos_1  & = & 0, \label{QuantAugI} \\
    \qrestr \prol    & = & \id_{-1}, \label{QuantAugII}
\end{eqnarray}
where we have written
$\id_{-1}:=\id_{\Grass \mathfrak g^*\otimes C^\infty(C)[[\lambda]]}$.
Moreover the map $\prol\qrestr$ is a projection onto
$\Grass \mathfrak g^*\otimes \mathcal F_C[[\lambda]]$ whose kernel is given by 
$\Grass \mathfrak g^*\otimes \qideal_C$.
In particular, we have the decomposition
\begin{equation}\label{QDirectSum}
    \Grass \mathfrak g^*\otimes C^{\infty}(M)[[\lambda]]
    = (\Grass \mathfrak g^*\otimes{\mathcal F}_{C})[[\lambda]]
      \oplus (\Grass \mathfrak g^*\otimes\qideal_{C}).
\end{equation}
\end{proposition}
\begin{proof}
Since $\qkos$ is a deformation of $\kos$ the formal series of
differential operators $B_1$ is well-defined. Moreover, we have
$\id_{0} = \prol\iota^* + \kos_1h_0 
= \prol\iota^* + \qkos_1h_0 +\lambda B_1h_0$
whence
\begin{equation} \label{QRestI}
    \id_0 = \prol\qrestr
            + \qkos_1 h_0(\id_0-\lambda B_1h_0)^{-1}.
\end{equation}
Because $h_0\prol = 0$ (Eq.~(\ref{WichtigeGleichung}))
we immediately see that 
$\qrestr\prol=\id_{-1}$ (since $\iota^*\prol=\id_{-1}$) which implies
that $\prol\qrestr$ is a projection. This fact together with
(\ref{QRestI}) entails that the product 
$\prol\qrestr\qkos_1 h_0(\id_0-\lambda B_1h_0)^{-1}$ vanishes, and
since $\prol$ is injective and $(\id_0-\lambda B_1h_0)^{-1}$ is
invertible we get $\qrestr\qkos_1h_0=0$. Multiplying $\qkos_1h_0$ by
$\kos_1$ from the right we obtain the following, using Lemma
\ref{KosLem} and writing $\id_1$ for 
$\id_{\Grass \mathfrak g^*\wedge \mathfrak g\otimes C^\infty(M)[[\lambda]]}$
\begin{equation}\label{QKosICrack}
    \qkos_1 h_0 \kos_1 = \qkos_1(\id_1-\kos_2h_1)
                       =   \qkos_1(\id_1-(\kos_2-\qkos_2)h_1),
\end{equation}
thanks to $0=\qkos_1\qkos_2$. But the difference $\kos_2-\qkos_2$ is a
multiple of $\lambda$ whence the right factor of the second equation
above is invertible which immediately implies $\qrestr\qkos_1 = 0$. 
The same argument (\ref{QKosICrack}) also shows that the image of
$\qkos_1$ is equal to the image of the projector 
$\qkos_1h_0(\id_{0}-\lambda B_1h_0)^{-1}$. It follows that
$\Grass\mathfrak g^*\otimes\qideal_C$ is the kernel of the
projection $\prol\qrestr$. Moreover, let $\phi\in C^\infty(C)[[\lambda]]$.
Then $\prol\qrestr\prol(\phi)=\prol(\phi)$ which shows that
the image of $\prol\qrestr$ is given by $\mathcal F_C[[\lambda]]$.
This proves the direct sum decomposition after tensoring with
the Gra{\ss}mann algebra.
Now suppose that there were another such deformed restriction
map, $\qrestr'$, satisfying the conditions of the Proposition.
Then (\ref{QuantAugII}) implies that $\prol\qrestr'$ is a
projection whose image is 
$\Grass\mathfrak g^*\otimes\mathcal F_C[[\lambda]]$.
Moreover, (\ref{QuantAugI}) entails that the image of 
$\qkos_1$, $\qideal_C$,
is contained in the kernel of this projection. Thanks to
the already proven direct sum decomposition (\ref{QDirectSum}) the
kernel of the projection $\prol\qrestr'$ decomposes into the direct
sum of the Gra{\ss}mann algebra tensor the quantum vanishing ideal and
the intersection of this kernel with 
$\Grass\mathfrak g^*\otimes \mathcal F_C[[\lambda]]$: this is
impossible by (\ref{QuantAugII}). Hence $\prol(\qrestr-\qrestr')=0$
which implies equality of the two quantum restriction maps by
injectivity of $\prol$. 
\end{proof}

We shall now define the augmented quantum BRST complex and the
deformed chain homotopies in a manner analogous to the classical case: 
let $\hat{\mathcal A}[[\lambda]]$ denote the \emph{augmented BRST complex}
$(\Grass \mathfrak g^*\otimes C^\infty(C)[[\lambda]])
\oplus \mathcal A[[\lambda]]$. Let 
$\hat{\qkos}:\hat{\mathcal A}[[\lambda]] 
\to \hat{\mathcal A}[[\lambda]]$
be equal to the quantised Koszul operator $\qkos$ on the quantum BRST
complex $\mathcal A^{\bullet,>0}[[\lambda]]$, equal to $\qrestr$ on 
$\mathcal A^{\bullet,0}[[\lambda]]$, and zero otherwise. The deformed
chain homotopies on the augmented quantised BRST complex will be
defined in the following 
\begin{proposition} \label{QHomotopieDef}
Define the quantised chain homotopy
$\hat{\qhomot}: \hat{\mathcal A}[[\lambda]] 
\to \hat{\mathcal A}[[\lambda]]$ 
in the following way: on the subspace 
$\Grass\mathfrak g^*\otimes C^\infty(C)[[\lambda]]$ set it equal to
the classical prolongation, i.e. $\hat{\qhomot}_{-1}:=\prol$ and on
the quantum BRST complex $\mathcal A[[\lambda]]$ set
$\hat{\qhomot}:=\qhomot$ where the latter is given as 
\begin{eqnarray}
    \qhomot_0 
    & := & h_0(\prol\qrestr+\qkos_1 h_0)^{-1} 
           = h_0(\id_0-\lambda B_1h_0)^{-1}, \label{QuantHomotNull} \\
    \qhomot_i 
    & := & h_i(h_{i-1}\qkos_i + \qkos_{i+1}h_i)^{-1}
           \qquad \forall i\geq 1, \label{QuantHomotSonst}
   \end{eqnarray}
where $B_1$ is as in Proposition~\ref{QAug}.
Furthermore we get the chain homotopy equation on the augmented
quantum BRST complex  
\begin{equation} \label{QuantHomotEq}
    \hat{\qkos} \hat{\qhomot} + \hat{\qhomot} \hat{\qkos} = \id,
\end{equation}
whence the quantised augmented Koszul complex also has trivial homology. 
\end{proposition}
\begin{proof}
Firstly, note that the fact that $\hat{\qkos}$ is a deformation of
$\hat{\kos}$ implies that $\hat{\qhomot}$ is a well-defined
deformation of the classical augmented chain homotopy
$\hat{h}$. Since $\hat{\qkos}^2=0$ it follows that 
$(\hat{h} \hat{\qkos} + \hat{\qkos} \hat{h})$ commutes with
$\hat{\qkos}$ which immediately implies
(\ref{QuantHomotEq}). Moreover, since the restriction of
$(\hat{h}\hat{\qkos}+\hat{\qkos} \hat{h})$ to $C^\infty(M)[[\lambda]]$
is equal to
\[
    \prol \qrestr + \qkos_1 h_0 = \id_0 -\lambda (\id_0-\prol\qrestr)B_1h_0,
\]
we see that the terms containing the projection $\prol\qrestr$ in
(\ref{QuantHomotNull}) vanish thanks to $h_0\prol = 0$, see 
(\ref{WichtigeGleichung}), 
which shows the second equation in (\ref{QuantHomotNull}). 
\end{proof}

The following lemma is a key to questions of bidifferentiability of the
reduced star product:
\begin{lemma} \label{DiffQRestr}
Using the geometric chain homotopy $h_0$ of Lemma \ref{KosLem} there
is a formal series of differential linear operators of $C^\infty(M)$,
$S:=\id_{C^\infty(M)}+\sum_{r=1}^\infty \lambda^r S_r$ where $S_r$
vanishes on constants for $r \ge 1$, such that
\[
    \qrestr = \iota^*\circ S.
\]
In particular, if $*$ is of Vey type then the order of each $S_r$ is
at most $r$. Moreover, under the additional assumptions of
Lemma~\ref{L31inv} $S$ can be chosen to be $G$-invariant.
\end{lemma}
\begin{proof}
It suffices to show that for any differential operator 
$D:\mathfrak g \otimes C^\infty(M)\rightarrow C^\infty(M)$ of order
$k$ there is a differential operator 
$D':C^\infty(M) \to C^\infty(M)$ of order $k+1$ such that the following
holds in the tubular neighbourhood of Lemma \ref{KosLem}:
\[
    \iota^*\circ D\circ h_0 = \iota^*\circ D'.
\]
But this is clear: the tubular neighbourhood is diffeomorphic to an open
subset of $C\times \mathfrak g^*$, hence using the classical momentum map $J$
as a global coordinate we see that $D$ takes the form
\[
    D = \sum_{s=0}^k \sum_{j,j_1,\ldots,j_s=1}^{\dim {\mathfrak g}}
        i(e^j){D^{(s)}}_{j j_1\cdots j_s}
        \frac{\partial^s}{\partial J_{j_1}\cdots \partial J_{j_s}},
\]
where ${D^{(s)}}_{jj_1\cdots j_s}$ are differential operators
$C^\infty(C)\rightarrow C^\infty(C)$ of
order $k-s$ which are smoothly parametrised by
$J_1, \ldots, J_n$, and where $e^1, \ldots, e^n$ is a basis of
$\mathfrak g^*$. Using formula (\ref{HomotNearC}) for $h_0$ near
$C$ we see that
\[
    D' = \sum_{s=0}^k \sum_{j,j_1,\ldots,j_s=1}^{\dim {\mathfrak g}}
         \frac{1}{s+1} {D^{(s)}}_{j j_1\cdots j_s}
         \frac{\partial^{s+1}}
         {\partial J_j\partial J_{j_1}\cdots \partial J_{j_s}}
\]
will satisfy the above equation and is clearly $G$-invariant if $D$ is 
$G$-equivariant. By induction, the composition of the
restriction $\iota^*$ and finitely many powers of $\lambda B_1h_0$
will be of the desired differential operator form with the correct
bounds for the orders. This will define the operators $S_r$, 
$r \geq 1$, of the asserted formula in the ($G$-invariant) tubular
neighbourhood. Since by the very definition of a star product $B_1$
vanishes on the constant functions, the $S_r$ will also 
vanish on the constants for $r\geq 1$. Multiplying the $S_r$, 
$r \geq 1$, by a suitable ($G$-invariant) bump function with support
in the tubular neighbourhood $S$ (and equal to $1$ in an open
neighbourhood of $C$) will give us a globally defined operator series
$S$ still satisfying the asserted equation. 
\end{proof}

We shall now need the quantum analogue of the Lie algebra
representation $\ClassLieC$ on the constraint surface to construct the
quantum analogue of the Chevalley-Eilenberg differential $\ce^c$. The
motivating classical equation is the identity (\ref{ClassConstrRep}),
and we set for all $\xi \in \mathfrak g$ 
\begin{equation}
\label{QuantConstrRep}
    \QuantLieC(\xi) := \qrestr\QuantLieM(\xi)\prol.
\end{equation}
where the representation $\QuantLieM$ is defined in (\ref{QuantLieMDef}).
\begin{lemma}
The map $\QuantLieC$ defines a Lie algebra representation of
$\mathfrak g$ on the space $C^\infty(C)[[\lambda]]$. Moreover, the
quantised restriction map $\qrestr$ induces a $\mathfrak g$-module
isomorphism of the $\mathfrak g$-module
$C^\infty(M)[[\lambda]]/\qideal_C$ (compare Lemma \ref{QIdealRep})
onto the $\mathfrak g$-module $C^\infty(C)[[\lambda]]$. More
precisely, we have the following identity for all $\xi\in\mathfrak g$: 
\begin{equation}\label{QIntertw}
    \QuantLieC(\xi) \qrestr = \qrestr \QuantLieM(\xi).
\end{equation}
\end{lemma}
\begin{proof}
Since the map $\prol\qrestr$ is a projection (see Proposition \ref{QAug}) whose
kernel (after restriction to ghostnumber zero) is equal to the quantum vanishing
ideal it follows from Lemma \ref{QIdealRep} that for all $\xi\in\mathfrak g$
\[
    \prol\qrestr \QuantLieM(\xi) \prol\qrestr = \prol\qrestr \QuantLieM (\xi).
\]
Hence for all $\xi,\eta\in\mathfrak g$
\[
    \prol \QuantLieC(\xi)\QuantLieC(\eta) \qrestr
    = \prol\qrestr \QuantLieM(\xi) \prol\qrestr \QuantLieM(\eta) 
      \prol\qrestr
    = \prol\qrestr \QuantLieM(\xi)\QuantLieM(\eta) \prol\qrestr,
\]
whence the representation identity follows since $\prol$ is injective
and $\qrestr$ is surjective. The rest of the lemma is clear thanks
to Proposition \ref{QAug}. 
\end{proof}

As in the classical case (\ref{CCEAnticom}) there is the following simple 
consequence for the corresponding quantum Chevalley-Eilenberg
differential, $\qce^c$, (see (\ref{CEDef}) for a definition) on
the constraint surface: 
\begin{equation}\label{QCEAnticom}
    \qce^c \qrestr = - \qrestr \qce.
\end{equation}
In complete analogy to the classical case we denote by
$\hat{\sbrso}$ the augmented quantum $BRST$ operator
\begin{equation}
\label{AugQuantOp}
    \hat{\sbrso} := \qce^c + 2\qrestr + \sbrso,
\end{equation}
where all the maps are defined to be zero on the domains on which they were
previously not defined. Clearly
\begin{equation}
\label{AugQuantEq}
    \hat{\sbrso}^2 = 0.
\end{equation}
The augmented quantum BRST complex is depicted in Figure \ref{QuantKomplex}.
%
%
%
%
\begin{figure}
\begin{center}
\qquad
{
\unitlength=0.86pt
\footnotesize
\begin{picture}(455.00,290.00)(0.00,0.00)
\put(0.00,20.00){\makebox(0.00,0.00){$\pi^*$}}
\put(180.00,280.00){\vector(-1,0){30.00}}
\put(410.00,50.00){\makebox(0.00,0.00){$\qkos_n$}}
\put(410.00,90.00){\makebox(0.00,0.00){$\qkos_n$}}
\put(410.00,130.00){\makebox(0.00,0.00){$\qkos_n$}}
\put(410.00,210.00){\makebox(0.00,0.00){$\qkos_n$}}
\put(410.00,250.00){\makebox(0.00,0.00){$\qkos_n$}}
\put(410.00,290.00){\makebox(0.00,0.00){$\qkos_n$}}
\put(445.00,260.00){\makebox(0.00,0.00){$\qce$}}
\put(445.00,220.00){\makebox(0.00,0.00){$\qce$}}
\put(445.00,100.00){\makebox(0.00,0.00){$\qce$}}
\put(445.00,60.00){\makebox(0.00,0.00){$\qce$}}
\put(350.00,60.00){\makebox(0.00,0.00){$\qce$}}
\put(350.00,100.00){\makebox(0.00,0.00){$\qce$}}
\put(350.00,220.00){\makebox(0.00,0.00){$\qce$}}
\put(350.00,260.00){\makebox(0.00,0.00){$\qce$}}
\put(455.00,80.00){\makebox(0.00,0.00){$\mathcal A^{1,n}[[\lambda]]$}}
\put(455.00,120.00){\makebox(0.00,0.00){$\mathcal A^{2,n}[[\lambda]]$}}
\put(455.00,200.00){\makebox(0.00,0.00){$\mathcal A^{n-2,n}[[\lambda]]$}}
\put(455.00,240.00){\makebox(0.00,0.00){$\mathcal A^{n-1,n}[[\lambda]]$}}
\put(455.00,280.00){\makebox(0.00,0.00){$\mathcal A^{n,n}[[\lambda]]$}}
\put(360.00,280.00){\makebox(0.00,0.00){$\mathcal A^{n,n-1}[[\lambda]]$}}
\put(360.00,240.00){\makebox(0.00,0.00){$\mathcal A^{n-1,n-1}[[\lambda]]$}}
\put(360.00,200.00){\makebox(0.00,0.00){$\mathcal A^{n-2,n-1}[[\lambda]]$}}
\put(360.00,120.00){\makebox(0.00,0.00){$\mathcal A^{2,n-1}[[\lambda]]$}}
\put(360.00,80.00){\makebox(0.00,0.00){$\mathcal A^{1,n-1}[[\lambda]]$}}
\put(455.00,40.00){\makebox(0.00,0.00){$\mathcal A^{0,n}[[\lambda]]$}}
\put(360.00,40.00){\makebox(0.00,0.00){$\mathcal A^{0,n-1}[[\lambda]]$}}
\put(285.00,160.00){\makebox(0.00,0.00){.....}}
\put(410.00,160.00){\makebox(0.00,0.00){.....}}
\put(455.00,90.00){\vector(0,1){20.00}}
\put(455.00,50.00){\vector(0,1){20.00}}
\put(455.00,250.00){\vector(0,1){20.00}}
\put(455.00,210.00){\vector(0,1){20.00}}
\put(455.00,170.00){\vector(0,1){20.00}}
\put(455.00,130.00){\vector(0,1){20.00}}
\put(360.00,130.00){\vector(0,1){20.00}}
\put(360.00,170.00){\vector(0,1){20.00}}
\put(360.00,50.00){\vector(0,1){20.00}}
\put(360.00,90.00){\vector(0,1){20.00}}
\put(360.00,210.00){\vector(0,1){20.00}}
\put(360.00,250.00){\vector(0,1){20.00}}
\put(425.00,40.00){\vector(-1,0){30.00}}
\put(425.00,80.00){\vector(-1,0){30.00}}
\put(425.00,120.00){\vector(-1,0){30.00}}
\put(425.00,200.00){\vector(-1,0){30.00}}
\put(425.00,240.00){\vector(-1,0){30.00}}
\put(425.00,280.00){\vector(-1,0){30.00}}
\put(330.00,40.00){\vector(-1,0){30.00}}
\put(330.00,80.00){\vector(-1,0){30.00}}
\put(330.00,120.00){\vector(-1,0){30.00}}
\put(330.00,200.00){\vector(-1,0){30.00}}
\put(330.00,240.00){\vector(-1,0){30.00}}
\put(330.00,280.00){\vector(-1,0){30.00}}
\put(10.00,0.00){\makebox(0.00,0.00){$C^{\infty}(M_{\mbox{\rm{\tiny {red}}}})[[\lambda]]$}}
\put(210.00,120.00){\makebox(0.00,0.00){$\mathcal A^{2,1}[[\lambda]]$}}
\put(210.00,80.00){\makebox(0.00,0.00){$\mathcal A^{1,1}[[\lambda]]$}}
\put(210.00,40.00){\makebox(0.00,0.00){$\mathcal A^{0,1}[[\lambda]]$}}
\put(120.00,120.00){\makebox(0.00,0.00){$\mathcal A^{2,0}[[\lambda]]$}}
\put(120.00,80.00){\makebox(0.00,0.00){$\mathcal A^{1,0}[[\lambda]]$}}
\put(120.00,40.00){\makebox(0.00,0.00){$\mathcal A^{0,0}[[\lambda]]$}}
\put(285.00,120.00){\makebox(0.00,0.00){.....}}
\put(285.00,80.00){\makebox(0.00,0.00){.....}}
\put(285.00,40.00){\makebox(0.00,0.00){.....}}
\put(110.00,100.00){\makebox(0.00,0.00){$\qce$ }}
\put(110.00,60.00){\makebox(0.00,0.00){$\qce$ }}
\put(200.00,100.00){\makebox(0.00,0.00){$\qce$ }}
\put(200.00,60.00){\makebox(0.00,0.00){$\qce$ }}
\put(0.00,120.00){\makebox(0.00,0.00){$\bigwedge^2\mathfrak{g}^*\otimes C^{\infty}(C)[[\lambda]]$}}
\put(10.00,80.00){\makebox(0.00,0.00){$\mathfrak{g}^*\otimes C^{\infty}(C)[[\lambda]]$}}
\put(10.00,40.00){\makebox(0.00,0.00){$C^{\infty}(C)[[\lambda]]$}}
\put(0.00,100.00){\makebox(0.00,0.00){$\qce^c$ }}
\put(0.00,60.00){\makebox(0.00,0.00){$\qce^c$ }}
\put(80.00,130.00){\makebox(0.00,0.00){\boldmath{$\iota^*$} }}
\put(80.00,90.00){\makebox(0.00,0.00){\boldmath{$\iota^*$} }}
\put(80.00,50.00){\makebox(0.00,0.00){\boldmath{$\iota^*$} }}
\put(165.00,50.00){\makebox(0.00,0.00){$\qkos_1$}}
\put(165.00,90.00){\makebox(0.00,0.00){$\qkos_1$}}
\put(165.00,130.00){\makebox(0.00,0.00){$\qkos_1$}}
\put(255.00,130.00){\makebox(0.00,0.00){$\qkos_2$}}
\put(255.00,90.00){\makebox(0.00,0.00){$\qkos_2$}}
\put(255.00,50.00){\makebox(0.00,0.00){$\qkos_2$}}
\put(10.00,10.00){\vector(0,1){20.00}}
\put(90.00,40.00){\vector(-1,0){30.00}}
\put(180.00,40.00){\vector(-1,0){30.00}}
\put(270.00,40.00){\vector(-1,0){30.00}}
\put(210.00,50.00){\vector(0,1){20.00}}
\put(120.00,50.00){\vector(0,1){20.00}}
\put(10.00,50.00){\vector(0,1){20.00}}
\put(270.00,80.00){\vector(-1,0){30.00}}
\put(180.00,80.00){\vector(-1,0){30.00}}
\put(90.00,80.00){\vector(-1,0){30.00}}
\put(10.00,90.00){\vector(0,1){20.00}}
\put(210.00,90.00){\vector(0,1){20.00}}
\put(120.00,90.00){\vector(0,1){20.00}}
\put(270.00,120.00){\vector(-1,0){30.00}}
\put(180.00,120.00){\vector(-1,0){30.00}}
\put(90.00,120.00){\vector(-1,0){30.00}}
\put(210.00,130.00){\vector(0,1){20.00}}
\put(120.00,130.00){\vector(0,1){20.00}}
\put(10.00,130.00){\vector(0,1){20.00}}
\put(255.00,290.00){\makebox(0.00,0.00){$\qkos_2$}}
\put(0.00,220.00){\makebox(0.00,0.00){$\qce^c$ }}
\put(0.00,260.00){\makebox(0.00,0.00){$\qce^c$ }}
\put(110.00,260.00){\makebox(0.00,0.00){$\qce$ }}
\put(110.00,220.00){\makebox(0.00,0.00){$\qce$ }}
\put(200.00,220.00){\makebox(0.00,0.00){$\qce$ }}
\put(200.00,260.00){\makebox(0.00,0.00){$\qce$ }}
\put(80.00,160.00){\makebox(0.00,0.00){.....}}
\put(170.00,160.00){\makebox(0.00,0.00){.....}}
\put(285.00,240.00){\makebox(0.00,0.00){.....}}
\put(285.00,200.00){\makebox(0.00,0.00){.....}}
\put(285.00,280.00){\makebox(0.00,0.00){.....}}
\put(270.00,200.00){\vector(-1,0){30.00}}
\put(270.00,240.00){\vector(-1,0){30.00}}
\put(270.00,280.00){\vector(-1,0){30.00}}
\put(120.00,170.00){\vector(0,1){20.00}}
\put(10.00,170.00){\vector(0,1){20.00}}
\put(210.00,170.00){\vector(0,1){20.00}}
\put(80.00,250.00){\makebox(0.00,0.00){\boldmath{$\iota^*$} }}
\put(80.00,210.00){\makebox(0.00,0.00){\boldmath{$\iota^*$} }}
\put(80.00,290.00){\makebox(0.00,0.00){\boldmath{$\iota^*$}}}
\put(210.00,210.00){\vector(0,1){20.00}}
\put(210.00,250.00){\vector(0,1){20.00}}
\put(210.00,200.00){\makebox(0.00,0.00){$\mathcal A^{n-2,1}[[\lambda]]$}}
\put(0.00,200.00){\makebox(0.00,0.00){$\bigwedge^{n-2}\mathfrak{g}^*\otimes C^{\infty}(C)[[\lambda]]$}}
\put(120.00,200.00){\makebox(0.00,0.00){$\mathcal A^{n-2,0}[[\lambda]]$}}
\put(0.00,240.00){\makebox(0.00,0.00){$\bigwedge^{n-1}\mathfrak{g}^*\otimes C^{\infty}(C)[[\lambda]]$}}
\put(0.00,280.00){\makebox(0.00,0.00){$\bigwedge^n\mathfrak{g}^*\otimes C^{\infty}(C)[[\lambda]]$}}
\put(165.00,250.00){\makebox(0.00,0.00){$\qkos_1$}}
\put(165.00,210.00){\makebox(0.00,0.00){$\qkos_1$}}
\put(255.00,210.00){\makebox(0.00,0.00){$\qkos_2$}}
\put(165.00,290.00){\makebox(0.00,0.00){$\qkos_1$}}
\put(255.00,250.00){\makebox(0.00,0.00){$\qkos_2$}}
\put(210.00,240.00){\makebox(0.00,0.00){$\mathcal A^{n-1,1}[[\lambda]]$}}
\put(120.00,240.00){\makebox(0.00,0.00){$\mathcal A^{n-1,0}[[\lambda]]$}}
\put(210.00,280.00){\makebox(0.00,0.00){$\mathcal A^{n,1}[[\lambda]]$}}
\put(120.00,280.00){\makebox(0.00,0.00){$\mathcal A^{n,0}[[\lambda]]$}}
\put(90.00,280.00){\vector(-1,0){30.00}}
\put(10.00,210.00){\vector(0,1){20.00}}
\put(90.00,200.00){\vector(-1,0){30.00}}
\put(120.00,210.00){\vector(0,1){20.00}}
\put(180.00,200.00){\vector(-1,0){30.00}}
\put(180.00,240.00){\vector(-1,0){30.00}}
\put(10.00,250.00){\vector(0,1){20.00}}
\put(90.00,240.00){\vector(-1,0){30.00}}
\put(120.00,250.00){\vector(0,1){20.00}}
\end{picture}}
\caption{The augmented quantum BRST complex.} \label{QuantKomplex}
\end{center}
\end{figure}
We keep the
notation $\hat{\qhomot}$ for $\prol+\qhomot$ on 
$\hat{\mathcal A}$. Moreover, let $\hat{\qce}$ denote the augmented
Chevalley-Eilenberg operator $\qce^c+\qce$.

We are now going to compute the cohomology of the BRST complex in
terms of the quantum Chevalley-Eilenberg cohomology on the constraint
surface: 
\begin{theorem}
\label{Qcohomcrack}
With the above notations and definitions we have the following:
\begin{enumerate}
\item The following map
      \begin{equation} 
      \label{QuantModHomDef}
          \hat{\qhomot}' := 
          \frac{1}{2}\hat{\qhomot} \big(\id +
          \frac{1}{2}(\hat{\qce}\hat{\qhomot} +
          \hat{\qhomot}\hat{\qce})\big)^{-1} 
      \end{equation}
      is a homotopy for the differential $\hat{\sbrso}$, i.e. 
      $\hat{\sbrso}\hat{\qhomot}'+\hat{\qhomot}'\hat{\sbrso} = \id$.
      The $\mathbb C[[\lambda]]$-linear map
      \begin{equation} \label{QuantKosIso}
          \QPsi: \qcohom^{(\bullet)}_\brsind (\mathcal A[[\lambda]])
          \to
          \qcohom^{\bullet}_\ceind(\mathfrak g, C^\infty(C)[[\lambda]]):
          \quad
          [a] \mapsto [\qrestr a]
      \end{equation}
      is an isomorphism with the following inverse (where 
      $[c] \in \qcohom^{\bullet}_\ceind(\mathfrak g,C^\infty(C)[[\lambda]])$): 
      \begin{equation}\label{QuantKosInv}
          \QPsi^{-1}: [c] \mapsto \left[2\hat{\qhomot}'c\right]
          = \sum_{k=0}^{n} \left(-\frac{1}{2}\right)^k 
            \left[(\qhomot\qce)^k\prol c\right].
      \end{equation}
\item The isomorphism $\QPsi$ turns the quantum Chevalley-Eilenberg
      cohomology on the constraint surface, 
      $\qcohom^{\bullet}_\ceind(\mathfrak g, C^\infty(C)[[\lambda]])$,
      into a $\mathbb Z$-graded associative algebra with unit. The
      multiplication (also called star product) $*$ for 
      $c_1, c_2 \in \Grass \mathfrak g^* \otimes C^\infty (C)[[\lambda]]$
      takes the following form:
      \begin{equation} \label{StaProQCCE}
          [c_1] * [c_2] := 
          4\left[\qrestr \left(
          (\hat{\qhomot}'c_1) \spr (\hat{\qhomot}'c_2)\right)\right]
      \end{equation}
\item The homomorphism $L:\qcohom_\brsind^0(\mathcal A[[\lambda]])\rightarrow 
      \qbideal_C/\qideal_C$ in Proposition \ref{BRSTDirac} is an
      isomorphism of associative algebras.
\item For 
      $c_1, c_2 \in \qcohom_\brsind^0 (\mathcal A[[\lambda]]) 
      = \left(C^\infty (C)[[\lambda]]\right)^{\mathfrak g}$ 
      formula (\ref{StaProQCCE}) simplifies to
      \begin{equation}
      \label{ReducedStarProd}
          c_1 * c_2 = 
          \qrestr \left( (\prol c_1) \spr (\prol c_2) \right).
      \end{equation}   
\end{enumerate}
\end{theorem}
\begin{proof}
The proof is entirely analogous to the proof of the classical
Proposition \ref{cohomcrack} since only purely cohomological
statements are needed. For the last two parts note that $\prol c$ is in
$\qbideal_C$ if $\qce^c c=0$ and use the proof of 
Proposition~\ref{BRSTDirac}. 
\end{proof}

Note that, as also in the classical case, one has by pure diagram
chase and the acyclicity of the quantum Koszul complex
$\qcohom_{\ceind}^\bullet (\qcohom_0^\kosind (\mathcal A[[\lambda]]))
\cong \qcohom^{(\bullet)}_\brsind (\mathcal A[[\lambda]])$, 
but we do not need this fact.

\vspace{0.5cm}

We have thus shown that the quantum BRST cohomology may be computed in terms
of the quantum Chevalley-Eilenberg cohomology on the constraint
surface. The natural question arises whether this latter cohomology is
related to a deformation of functions on the classical reduced phase,
at least at level ghost number zero. In the next section we shall see
in a simple counter example that this will in general \emph{not} be the
case: there are in general `fewer' quantum invariants than classical
invariants. But there will be a large class of positive situations in
Section \ref{ProGroupSec}. 
In order to give a precise notion of `fewer invariants' we shall need
the following general definition of a 
\emph{deformation of a subspace}: 
\begin{definition}
Let $E: W \hookrightarrow V$ be a subspace of the vector space $V$
over some field $k$. A $k[[\lambda]]$-submodule $\qW$ of
$V[[\lambda]]$ is called a deformation of $W$ if there exists a
deformation of the canonical embedding 
$E$ i.e.~a formal power series of linear maps 
\begin{equation}
    \qE = \qE_0 + \sum_{r=1}^\infty \lambda^r \qE_r: 
    W[[\lambda]] \to V[[\lambda]]
\end{equation}
with $\qE_0 = E$, such that $\qE(W[[\lambda]]) = \qW$. 
\end{definition}
Note that it is not hard to see that in algebraic language the above deformations 
exactly correspond to those $\lambda$-adically closed primary submodules of 
$V[[\lambda]]$ whose associated prime ideal vanishes
(see e.g. \cite[p.~434]{Jac89} for a definition) or satisfy the equivalent and more
practical condition that $\lambda v\in \qW$
implies $v\in \qW$ for all $v\in V[[\lambda]]$, but we shall not need this in
the sequel.

We now give a precise definition
of a \emph{consistent quantum reduction}:
\begin{definition}
\label{ConsistentRedDef}
Let $(M, *, \mathfrak g, \qmm, C)$ be a Hamiltonian quantum 
$\mathfrak g$-space regular constraint surface and let 
$(M, \omega, \mathfrak g, J, C)$ be the corresponding classical
limit. We shall say that $(M, *, \mathfrak g, \qmm, C)$ admits a
consistent quantum reduction if and only if 
$\qcohom_\ceind^0 (\mathfrak g, C^\infty (C)[[\lambda]])$
is a deformation of
$C^\infty (C)^{\mathfrak g}[[\lambda]]$.
\end{definition}

Note that the quantum vanishing ideal $\qideal_C$ is a deformation of
the classical vanishing ideal defining $\qE$ by the restriction of
$\id - \prol \qrestr$ to the classical vanishing ideal (where
Eqn. (\ref{WichtigeGleichung}) is used). However, $\qbideal_C$ is not
necessarily a deformation of $\mathcal B_C$. 

%
%

\section{An example and counter example}
\label{ExampleSec}

This section is divided in three parts: first, we give an example of a 
classical phase space reduction (with proper group action and regular
values of the momentum) for $M := T^*S^1 \times T^* S^1$. Then we
shall show how our method for quantum reduction works for a strongly
invariant star product. Finally, we give an explicit example of a non
strongly invariant star product which cannot be reduced.

We consider the phase space $M = T^*S^1 \times T^* S^1$. A point in
$M$ is described by the quadruple $(z,p,w,J)$ where 
$z = e^{\im\varphi}$ and $w = e^{\im\psi}$ are the coordinates along 
$S^1$ and $p$ resp. $J$ are the corresponding momentum variables. The 
symplectic form is then given by 
$\omega = - d\varphi \wedge dp + d\psi \wedge dJ$ where $d\varphi$
resp. $d\psi$ are the non-exact global one-forms on the corresponding
circle. As group action we take the
$U(1)$-action induced by the coordinate function $J$ which thus serves 
as momentum map. The constraint surface
$C = J^{-1} (\{0\})$ is simply given by $T^*S^1 \times S^1$. On
$C$ we use the coordinates $(z,p,w)$, then $C$ is embedded in 
$M$ by $\iota: (z,p,w) \mapsto (z,p,w,0)$. Note that the group action
is clearly proper and free. Finally, the reduced phase space
$M_{\rm red}$ is given by $T^*S^1$ with coordinates $(z,p)$ and the
projection $\pi: C \to M_{\rm red}$ is just $(z,p,w) \mapsto (z,p)$. 
Since the Lie algebra $\mathfrak u (1)$ is one-dimensional we can
identify $\mathfrak u(1) \otimes C^\infty (M)$ with $C^\infty (M)$ and 
similar also 
$\mathfrak u^*(1) \wedge \mathfrak u(1) \otimes C^\infty (M)$ and 
$\mathfrak u^*(1) \otimes C^\infty (M)$. Then the classical Koszul
operator $\kos_1$ is just the multiplication with $J$, i.e.
\begin{equation}
\label{ClassKoszulEx}
    (\kos_1 f) (z,p,w,J) = J f(z,p,w,J)
\end{equation}
for $f \in C^\infty (M)$. As prolongation we simply use 
\begin{equation}
\label{ClassProlEx}
    (\prol u)(z,p,w,J) = u(z,p,w),
\end{equation}
where $u \in C^\infty (C)$. Due to the requirement
$\id = \prol \iota^* + \kos_1 h_0$ the chain homotopy $h_0$ is given
by the difference quotient
\begin{equation}
\label{ClassHEx}
    (h_0 f)(z,p,w,J) = \frac{f (z,p,w,J) - f (z,p,w,0)}{J}
\end{equation}
for $J \ne 0$, which is smoothly extended by the first partial
derivative in the $J$-direction on $C$. Then indeed
$\id = \prol \iota^* + \kos_1 h_0$. Finally the classical Lie
algebra action on $C^\infty (M)$ is given by the Poisson bracket with
$J$ which can be computed explicitly,
\begin{equation}
\label{ClassLieActionM}
    \ClassLieM f = \{J, f\} = - \frac{\partial}{\partial \psi} f,
\end{equation}
where $f \in C^\infty (M)$. The action on $C^\infty (C)$ is similarly
given by
\begin{equation}
\label{ClassLieActionC}
    \ClassLieC u = - \frac{\partial}{\partial \psi} u,
\end{equation}
where $u \in C^\infty (C)$. Obviously 
$\ClassLieC = \iota^* \ClassLieM \prol$.

Now let us describe the quantised version. Firstly we consider the
case of a strongly invariant star product. We take the following
star product (anti-standard ordered in the first variables and
standard ordered in the second variables) on $M$
\begin{equation}
\label{StandardStarProd}
    f * g = \mu \circ \eu^{\frac{\lambda}{\im} \left(
            \frac{\partial}{\partial \varphi} \otimes
            \frac{\partial}{\partial p} 
            +
            \frac{\partial}{\partial J} \otimes 
            \frac{\partial}{\partial \psi} \right)}
            f \otimes g
\end{equation}
(where $\mu (f \otimes g) = fg$ denotes the pointwise multiplication)
which is indeed an associative deformation of $\mu$ with the correct
Poisson bracket in first order. Now the quantised Koszul operator
$\qkos_1$ is given by right multiplication with $J$. Due to the
particular form of $*$ it is simply the undeformed multiplication
\begin{equation}
\label{qKoszulStd}
    \qkos_1 f = f * J = fJ = \kos_1 f.
\end{equation}
Thus we do not have to deform the chain homotopy $h_0$ and the
restriction $\iota^*$ and simply have 
$\id = \prol \iota^* + \qkos_1 h_0$. Finally we compute the quantum Lie
algebra action by taking $*$-commutators with $J$. It turns out that
$*$ is strongly invariant, i.e. we have
\begin{equation}
\label{StdLieActionM}
    \QuantLieM f = \frac{1}{\im\lambda} \ad_* (J) f 
                 = \frac{1}{\im\lambda} (J*f - f*J) 
                 = \{J,f\} 
                 = -\frac{\partial}{\partial \psi} f = \ClassLieM f,
\end{equation}
and similarly $\QuantLieC = \iota^* \QuantLieM \prol = \ClassLieC$ for
the quantum Lie algebra action on $C$. Thus the invariant functions on 
$C$ are clearly in bijection to the functions on $M_{\rm red}$ via
$\pi^*$ and hence this Hamiltonian quantum $\mathfrak u(1)$-space
allows a consistent quantum reduction by setting 
$\qE = \qE_0 = E$. The 
star product $\starred$ on the reduced phase space is well-defined by 
$\pi^*(u \starred v) = \iota^* ((\prol \pi^* u) * (\prol \pi^* v))$,
which explicitly yields the expected result
\begin{equation}
\label{StandardStarRed}
    u \starred v = \mu \circ \eu^{\frac{\lambda}{\im}   
                   \frac{\partial}{\partial \varphi} \otimes
                   \frac{\partial}{\partial p}} \,
                   u \otimes v.
\end{equation}

Next we shall consider what happens if one does not use a strongly
invariant star product. Note that for a one-dimensional Lie algebra
any star product is necessarily covariant. We shall even use a star
product $\tilde *$ equivalent to $*$. Consider the equivalence
transformation
\begin{equation}
\label{EquivEx}
    S = \exp \left(\lambda P \frac{\partial}{\partial J} \right)
\end{equation}
where $P \in C^\infty (M)$ is a smooth function of the momentum $p$
alone. Then define
\begin{equation}
\label{PervStarDef}
     f \pervstar g = S \left( S^{-1}f * S^{-1}g\right),
\end{equation}
and obtain $SJ = J - \lambda P$ since 
$\frac{\partial}{\partial J} P = 0$. A straightforward computation
yields that this star product is in general no longer strongly
invariant but
\begin{equation}
\label{PervLieActionM}
    \tilde\QuantLieM f 
    = \frac{1}{\im\lambda} (J \pervstar f - f \pervstar J)
    = \{J, f\} + \im \, \ad_* (P) f,
\end{equation}
where $\ad_*$ denotes the commutator with respect to $*$. Here one uses 
the fact that $S$ commutes with differentiation in the $\varphi$-
and $\psi$-direction. Note that $\ad_* (P)$ differentiates only in 
the $\varphi$-direction and has coefficients depending on $p$ only. 
The quantum Koszul operator for this star product can be computed also 
explicitly yielding
\begin{equation}
\label{qKoszulPerv}
    \tilde{\qkos}_1 f = f \pervstar J
    = \kos f + \lambda fP - \lambda f * P,
\end{equation}
where we used the fact that $S$ is an automorphism of the undeformed
product and acts trivially on functions which do not depend on $J$. 
The operator $B_1$ is then given by
\begin{equation}
\label{BPerv}
    B_1 f = f * P - fP = - \ad_* (P) f.
\end{equation}
Let us now compute the quantum restriction 
$\qrestr = \iota^* (\id - \lambda B_1 h_0)^{-1}$ more explicitly. 
First notice that $B_1$ commutes with $h_0$.
Then a simple computation
using the Taylor expansion in $J$ around $J=0$ yields
\begin{equation}
\label{RestrictHk}
    \iota^* \; h_0^k = \frac{1}{k!} \; \iota^* \; 
                       \frac{\partial^k}{\partial J^k}
\end{equation}
for all $k \in \mathbb N$. 
Together with the fact that $\frac{\partial}{\partial J}$ 
commutes with $B_1 = - \ad_*(P)$, too, we finally obtain
\begin{equation}
\label{QuantRestPerv}
    \qrestr = \iota^* 
    \exp\left(-\lambda\ad_*(P)\frac{\partial}{\partial J}\right),
\end{equation}
illustrating Lemma~\ref{DiffQRestr}.
This enables us to compute the quantum Lie algebra action on $C$
explicitly 
\begin{equation}
\label{PervLieActionC}
    \tilde\QuantLieC u = -\frac{\partial u}{\partial \psi} 
    + \im \, \ad_* (P) u,
\end{equation}
where we can use the operator $\ad_* (P)$ also for functions on $C$
since it only contains differentiation in the $\varphi$-direction and
multiplication by functions depending on $p$ due to our particular
choice of $P$. Now let us use a more particular function for $P$,
namely we take $P (z,p,w,J) = p$ then 
$\ad_*(P) = \im\lambda \frac{\partial}{\partial \varphi}$. 
Now we consider the functions on the constraint surface $C$ which are
invariant under this quantum action $\tilde \QuantLieC$. Let 
$u = \sum_{r=0}^\infty \lambda^r u_r \in C^\infty (C)[[\lambda]]$ then 
$\tilde\QuantLieC u = 0$ implies in lowest order 
$\frac{\partial u_0}{\partial \psi} = 0$ as classically expected. Thus 
$u_0$ does not depend on $\psi$, i.e.~is constant along the second
circle in $C = T^*S^1 \times S^1$. The next order then implies
\begin{equation}
\label{KnockOutCond}
    \frac{\partial u_1}{\partial \psi} 
    = - \frac{\partial u_0}{\partial \varphi} 
    \quad
    \textrm{ resp. }
    \quad
    \frac{\partial u_r}{\partial \psi} 
    = - \frac{\partial u_{r-1}}{\partial \varphi} 
\end{equation}
for $r \ge 1$. Since $S^1$ is compact we can integrate this equation 
over the second circle. Together with the fact that $u_0$ does not
depend on $\psi$ this yields that 
$\frac{\partial u_0}{\partial \varphi}$ has to vanish too. By
induction we conclude that $\tilde\QuantLieC u = 0$ if
and only if $u$ does only depend on $p$ but not on $\varphi$. Thus we
do \emph{not} obtain all functions on the reduced phase space in this
case but only those which are constant along the $\varphi$-direction.
Thus we conclude that in this example there is no consistent quantum
reduction. 

%
%

\section {Proper group actions and other nice cases}
\label{ProGroupSec}

In view of the examples given in the preceding section we should like to
present some classes of Hamiltonian quantum $\mathfrak g$-spaces
$(M, *, \mathfrak g, \qmm, C)$ with regular constraint surface
which allow a consistent quantum reduction. Moreover we shall specify
how the star product for the reduced phase space $M_{\rm red}$ 
looks whenever $M_{\rm red}$ exists as manifold. 

The first class of examples is the class of 
\emph{proper Hamiltonian $G$-spaces and strongly invariant star products}:
\begin{theorem}
Let $(M, \omega, G, J, C)$ be a Hamiltonian $G$-space with regular
constraint surface where the Lie group $G$ is connected and acts
properly on $M$. Moreover, pick any strongly invariant
star product $*$ on $M$ (whose existence is assured by Fedosov's
Theorem), and set $\qmm:=J$. Furthermore choose $G$-equivariant chain
homotopies $h$ and a $G$-equivariant prolongation $\prol$ (whose
existence is assured by Lemma \ref{L31inv}). Then we have the following:
\begin{enumerate}
\item The quantum representations $\QuantLieM$ and $\QuantLieC$ of the Lie
      algebra $\mathfrak g$ of $G$ are equal to the corresponding classical
      representations $\ClassLieM$ and $\ClassLieC$ which implies that the
      quantum and classical Chevalley-Eilenberg differentials are equal,
      i.e.~$\qce = \ce$ and $\qce^c = \ce^c$. Hence
      \[
          \qcohom^\bullet_\ceind(\mathfrak g, C^\infty(C)[[\lambda]])
          =
          H^\bullet_\ceind(\mathfrak g, C^\infty(C))[[\lambda]],
      \]
      whence $(M, *, \mathfrak g, \qmm, C)$ admits a consistent
      quantum reduction. Moreover, we have 
      $\hat{\qce}\hat{\qhomot}+\hat{\qhomot}\hat{\qce}=0.$

\item The star product on the Chevalley-Eilenberg cohomology
      $H_\ceind(\mathfrak g, C^\infty(C))[[\lambda]]$ is given by the
      following simplified formula for $c_1,c_2\in \Grass \mathfrak
      g^*\otimes C^\infty(C)[[\lambda]]$ with $\ce c_1 = 0 = \ce c_2$:
      \[
          [c_1]*[c_2] = 
          \left[\qrestr((\prol c_1)\spr (\prol c_2))\right].
      \]
\item Suppose in addition that the $G$-action on
      $C$ is free, and let
      $\pi:C\rightarrow M_{\rm red}:=C/G$ the canonical projection
      onto the reduced phase space $M_{\rm red}$. Then there is the
      following simple formula for the reduced star product $\starred$
      of two functions 
      $\phi_1,\phi_2\in C^\infty(M_{\rm red})[[\lambda]]$
      which is induced by the above construction:
      \begin{equation}
      \label{SimpleStarFormula}
          \pi^*(\phi_1 \starred \phi_2) 
          = \qrestr((\prol \pi^* \phi_1)*(\prol \pi^*\phi_2)).
      \end{equation}
      Suppose that the prolongation and the chain homotopies are
      geometric then the star product $\starred$ is bidifferential and
      if $*$ is even of Vey type then $\starred$ is of Vey type,
      too. 

      Moreover, for any two $G$-invariant functions $f_1$ and
      $f_2$ in $C^\infty (M)[[\lambda]]$ we have the following: let ${f_1}_{\rm red}$ and
      ${f_2}_{\rm red}$ be the unique functions of 
      $C^\infty(M_{\rm red})[[\lambda]]$ such that 
      $\qrestr f_k =: \pi^*{f_k}_{\rm red}$ for $k=1,2$. Then
      \begin{equation}
           \pi^*({f_1}_{\rm red}\starred {f_2}_{\rm red})
                      = \qrestr(f_1*f_2).       
      \end{equation}

\item The choice of a different $G$-invariant geometric prolongation
      and different geometric chain homotopies yields in general a
      different but equivalent reduced star product. 
\end{enumerate}
\end{theorem}
\begin{proof}
The first part is clear: one has equality of quantum and classical
representations
since $*$ is strongly invariant, since $\prol$ 
intertwines the classical $\mathfrak g$-actions on $M$ and on $C$ in
formula (\ref{QuantConstrRep}), and since $\qrestr\prol$ is the
identity on $C^\infty(C)[[\lambda]]$ by (\ref{QuantAugII}). Moreover
$\hat{\qce}$ anticommutes with $\hat{\qhomot}$ because the quantized
Koszul operator, the quantized restriction and the quantized chain 
homotopies are clearly $G$-equivariant since their classical 
counterparts are $G$-equivariant.
The second part is an immediate consequence of formula (\ref{StaProQCCE}) and
the fact that $\hat{\qhomot}'=\frac{1}{2}\qhomot$ by the above.
The formulas in the third part are simple consequences of part ii.)
and of the fact that 
$H^0_\ceind(\mathfrak g, C^\infty(C))=\pi^*C^\infty(M_{\rm red})$, see
the Propositions \ref{ReducClas} and \ref{ZeroCohomClass}. In order to
prove the differentiability properties of the reduced star product we
consider 
around some arbitrarily chosen point $x\in M_{\rm red}$ an open chart
$U$ such that the principal fibre bundle $C$ trivialises, hence
$\pi^{-1}(U)$ is $G$-equivariantly diffeomorphic to $U\times G$. Taking
the $G$-equivariant tubular neighbourhood of $C$ restricted to
$\pi^{-1}(U)$ we finally get a chart domain of $M$ diffeomorphic to
an open subset of $U\times G\times \mathfrak g^*$. In both cases
above the star product restricted to this domain is a series of
bidifferential operators. Since $\prol\pi^*\phi_1$ simply corresponds to
a function on this domain not depending on the coordinates on 
$G\times\mathfrak g^*$ and since the quantum restriction map $\qrestr$ is equal
to $\iota^*\circ S$ where $S$ is a series of differential
operators, or differential operators of order bounded by the order of
$\lambda$, respectively, according to Lemma \ref{DiffQRestr}, we
immediately see in co-ordinates that the resulting bilinear operators of the
reduced star products are bidifferential operators of the same order
which proves the two asserted properties at the same time. The last formula
follows from the fact that $f-\prol\pi^*f_{\rm red}$ is in the intersection
of the quantum vanishing ideal with the space of all smooth
complex-valued $G$-invariant functions, $C^\infty(M)^G[[\lambda]]$:
this intersection is a two-sided ideal of the subalgebra $C^\infty(M)^G[[\lambda]]$
which follows from the explicit form of the
elements of $\qideal_C$ (see Proposition \ref{QIdealRep}) and the fact that
$\langle J,\xi\rangle * f=f*\langle J,\xi\rangle $ for all 
$\xi \in \mathfrak g$ and all $f\in C^\infty(M)^G[[\lambda]]$ since $*$ is
strongly invariant. Fourthly, let $\prol'$ and $\qrestr'$ be another
choice and let $S$, $S'$ be $G$-invariant differential operators on $M$
according to Lemma~\ref{DiffQRestr} such that $\qrestr = \iota^* \circ S$ and
$\qrestr' = \iota^*\circ S'$. Then the linear map 
$\Phi: C^\infty (M_{\rm red})[[\lambda]] \to C^\infty (M_{\rm red})[[\lambda]]$ 
defined by 
$\Phi (\varphi):= (S^{-1} S' \prol \pi^* \varphi)_{\rm red}$ is easily
seen 
to be an algebra isomorphism of 
$(C^\infty (M_{\rm red})[[\lambda]],*_{\rm red})$ 
onto
$(C^\infty (M_{\rm red})[[\lambda]],*'_{\rm red})$ using the fact that 
$(\cdot)_{\rm red}$ is a homomorphism and that 
$f_{{\rm red}'}=(S^{-1} S'f)_{\rm red}=\Phi f_{\rm red}$ for each 
$f\in C^\infty(M)^G[[\lambda]]$. The fact that $\Phi$ is a formal
series of differential operators is shown by a local consideration
analogous to the one in part three.
\end{proof}

\begin{remark} 
The associativity of the explicit formula (\ref{SimpleStarFormula})
can also be seen more directly: note that the space of $G$-invariant
functions $C^\infty(M)^G[[\lambda]]$ on $M$, is a sub-algebra of
$(C^\infty(M)[[\lambda]],*)$ due to the $G$-invariance of $*$.
Furthermore, since the projection $\prol\qrestr$ is $G$-equivariant we
have the decomposition
\[
  C^\infty(M)^G[[\lambda]]
      =
       (C^\infty(M)^G[[\lambda]]\cap \mathcal F_C[[\lambda]])
                    \oplus
          (C^\infty(M)^G[[\lambda]]\cap \qideal_C).
\]
As has been shown above the space $C^\infty(M)^G[[\lambda]]\cap \qideal_C$ 
is a two-sided ideal in $C^\infty(M)^G[[\lambda]]$. 
Hence the $\mathcal F_C$-component of the star product
of two $G$-invariant functions yields an associative multiplication.
We use results on Koszul homology
to prove the directness of the above decomposition. This point of view
has been used by Schirmer to compute a 
star product on complex Gra{\ss}mann manifolds, see \cite{Schir97,Schir98}.
\end{remark}

The second class of examples is based on the following
\begin{theorem}
Let $(M, *, \mathfrak g, \qmm, C)$ be a Hamiltonian quantum
$\mathfrak g$-space with regular constraint surface.
Suppose in addition that the first classical Chevalley-Eilenberg
cohomology group on the constraint surface, 
$H^1_\ceind (\mathfrak g,C^\infty(C))$, vanishes. 
Then the Hamiltonian quantum
$\mathfrak g$-space allows a consistent quantum reduction.
\end{theorem}
\begin{proof}
Write $\qce^c=\sum_{r=0}^\infty \lambda^r \qce^c_r$ where $\qce^c_r$
are linear endomorphisms of $\Grass g^*\otimes C^\infty(C)$ and
$\qce^c_0=\ce^c$. We consider the equation
\[
    \qce^c \phi'=0
\]
where $\phi'=\sum_{r=0}^\infty \lambda^r \phi_r$ is in 
$C^\infty(C)[[\lambda]]$. Its solvability is a standard argument for
deformed cohomology operators (as for instance in the proof of the
existence of the Fedosov construction): Choose a vector subspace
$\mathcal C \subset C^\infty(C)$ such that 
$C^\infty(C)=C^\infty(C)^{\mathfrak g}\oplus \mathcal C$. 
We construct the maps 
$E_r, r\in\mathbb N: C^\infty(C)^{\mathfrak g} \to C^\infty(C)$ 
as follows: let $\phi\in C^\infty(C)^{\mathfrak g}$. Hence
$\ce^c \phi=0$ and the above equation is satisfied up to order zero.
Set $\qE_0$ equal to the canonical injection into $C^\infty(C)$ and
suppose that the maps $\qE_r$ have already been constructed up to
order $s \in \mathbb N$ such that the image of $\qE_r$ is contained in
the complementary space $\mathcal C$ for all $1\leq r\leq s$,
and such that $\phi^{(s)}:=\sum_{r=0}^s\lambda^r \qE_r(\phi)$ solves
the above equation up to order $s$, i.e. $(\qce^c\phi^{(s)})_r=0$ for
all $0\leq r\leq s$. Consider
$(\qce^c\phi^{(s)})_{s+1}$. Since obviously $\qce^c(\qce^c\phi^{(s)})=0$ we
get $0=(\qce^c(\qce^c\phi^{(s)}))_{s+1}=\ce^c(\qce^c\phi^{(s)})_{s+1}$.
But since $H^1_\ceind(\mathfrak g,C^\infty(C))$ vanishes there must be a
$\phi_{s+1}\in C^\infty(C)$ such that
\[
    (\qce^c\phi^{(s)})_{s+1}=-\ce^c\phi_{s+1}.
\]
Clearly $\phi_{s+1}$ is unique up to addition of an arbitrary element 
$\psi_{s+1}\in H^0_\ceind(\mathfrak g,C^\infty(C))$. Let
$\qE_{s+1}(\phi)$ be the unique such element $\phi_{s+1}$ in 
$\mathcal C$. This clearly defines a linear map 
$\qE_{s+1}: C^\infty(C)^{\mathfrak g} 
\to \mathcal C\subset C^\infty(C)$, 
and the above equation is clearly equivalent to the equation 
$0 = (\qce^c\phi^{(s+1)})_{s+1}$. Proceeding this
way by induction we construct all the maps $\qE_r$ such that 
$\phi' = \qE \phi = \sum_{r=0}^\infty\lambda^r \qE_r (\phi)$ solves 
$\qce^c\phi'=0$ for arbitrary 
$\phi\in C^\infty(C)^{\mathfrak g}$. Since on the other hand every
solution to $\qce^c\phi'=0$ can obviously constructed that way the
theorem is proved.
\end{proof}

The cohomological criterion in the preceding theorem can be made more
explicit in case the Hamiltonian quantum $\mathfrak g$-space comes from
a Hamiltonian $G$-space with existing reduced phase space:
\begin{lemma}
Let $(M,\omega,G,J,C)$ be a Hamiltonian $G$-space with regular
constraint surface where the connected Lie group freely and properly
acts on $C$. Then $H^1_\ceind(\mathfrak g,C^\infty(C))=0$ if and only
if the first de Rham cohomology group of the Lie group $G$ vanishes.
\end{lemma}
\begin{proof}
We know that $C$ is a principal fibre bundle over the reduced phase
space $M_{\rm red}:=C/G$ with structure group $G$. Observe that
$H^\bullet_{\mbox{\tiny de~Rham}}(G) 
\cong H^\bullet_\ceind(\mathfrak g,C^\infty(G))$
using right invariant differential forms as a global basis of all
differential forms on $G$ and letting $\mathfrak g$ act on
$C^\infty(G)$ by Lie derivatives of right invariant vector
fields. Moreover $\ce^c$ clearly commutes with multiplication by
pull-backs of smooth functions on the reduced phase space by means of
the projection $\pi$. Hence, using partitions of unity on the reduced
space it suffices to prove the assertion for cochains supported in
bundle charts diffeomorphic to $U\times G$. If 
$H^1_\ceind(\mathfrak g,C^\infty(C))=0$ choose a closed one-form
$\beta$ on $G$, i.e. an element in $\mathfrak g^*\otimes C^\infty(G)$,
multiply it with a bump function $b$ being equal to $1$ at some point
$x$ on the reduced space to obtain a $\ce^c$-closed element $\beta'$
in $\mathfrak g^*\otimes C^\infty(C)$.
Since $\beta'=\ce^c\phi$ then $\beta = d\phi(x,\cdot)$, showing one
implication. For the reverse implication note that every
$\beta'\in\mathfrak g^*\otimes C^\infty(C)$ which is $\ce^c$-closed gives
rise to a $d$-closed one-form $\beta'(u,\cdot)$ on $G$ smoothly
parametrised by $u\in U$. For instance by means of line integrals on
$G$ along arbitrary smooth paths starting at the neutral element we
can choose a smooth function $\phi \in C^\infty(U\times G)$ such that
$\beta'(u,\cdot) = d\phi(u,\cdot)$. Hence $\beta'$ is equal to 
$\ce^c \phi$, and this will establish the reverse implication.
\end{proof}

%
%

\section {Outlook and open problems}
\label{OutSec}

In this section we list some open questions arising with our approach
to BRST cohomology and give an outlook to future work we plan to do.

The first remark concerns the principal question of covariant star 
products. Obviously our approach relies crucially on the 
\emph{existence of an at least quantum covariant star product}
for the given Lie group or
Lie algebra action. Though in the literature many examples of the
existence of (quantum) covariant star products are known,
it seems that there is no general theorem on the existence of
covariant star products. Neither it is known whether there are
principal obstructions for the construction of covariant star products.
Thus it would be highly desirable to find here more concrete
conditions whether a covariant star product exists or even prove that
there are no obstructions at all beside the existence of a classical
momentum map. An answer to this question would have many useful
consequences not only for the BRST method. In particular the general
existence of a covariant star product would imply the 
\emph{quantum integrability} of any classically integrable system where 
one simply has a Hamiltonian $\mathbb R^n$-action.

Secondly, in view of the counter example in Section~\ref{ExampleSec},
there is the following open question: For any Hamiltonian quantum 
$\mathfrak g$-space with regular constraint surface 
$(M, *, \mathfrak g, \qmm, C)$ does there always exist an equivalent
star product $*'$ on $M$ such that $(M, *', \mathfrak g, \qmm, C)$ is
a Hamiltonian quantum $\mathfrak g$-space with regular constraint
surface allowing for a consistent quantum reduction?
Moreover, it may happen that equivalent but different star products on
$M$ induce non-equivalent star products on $M_{\rm red}$. Such effects
are known from phase space reduction of star products in several
examples, see e.g.~\cite{Wal98a,Gloe98a,Gloe98b}, and thus have to be
expected in the BRST framework as well. Finally it
would be interesting to see how one can compute the Deligne class
(see \cite{Del95} and \cite{GR99} for definitions)
of the reduced star product in terms of the BRST algebra. As a conversation
of M.B. with G.~Halbout and, independently, a proposal of the referee
suggests this will probably be related to a formulation of the
(BRST) reduction in terms of a suitable `equivariant cohomology'.

As third open problem we would like to mention that the formalism we
have presented here only deals with the 
\emph{observable algebra}. Though the observables are the primary
object of deformation quantization one has to deal with the question
how the \emph{states} of the considered physical system can be
described. In \cite{BW98a} the concept of
formally positive functionals and their GNS representations in
deformation quantization has been introduced and shown to be a
physically reasonable concept for the description of states. 
Thus we plan in a forthcoming paper to consider GNS
representations of the BRST algebra. Here some subtleties arise from
the fact that the $^*$-involution of complex conjugation is
now $\mathbb Z_2$-graded where we have to use the Weyl ordered case. 
Thus the concept of positive functionals and 
their GNS representation has to be modified. The main problem is to
find conditions for $\mathbb C[[\lambda]]$-linear functionals of the
BRST algebra such that firstly the GNS representation of the BRST
algebra induces indeed a representation of the algebra of smooth
functions on the reduced phase space, and, secondly, to guarantee the
positivity of the inner product of the representation space of the 
reduced observables. The problem of a positive inner product is known
in the usual approach to BRST quantization and the framework of
formal GNS representations seems to be a promising way to construct
positive inner products.

After the first version of this paper appeared on the internet
M.~Semenov-Tyan-Shanskii proposed to M.B. to investigate the more
general problem of a Lie-Poisson group action with a `noncommutative
momentum map' taking values in the dual group: here a (quantum) BRST
formulation of the phase space reduction does not yet seem to have
been achieved.

Fifthly, apart from BRST theory the more general problem of finding
\emph{deformed or quantised analogues of classical sub-manifolds} seems
to be interesting and also occurring (see also e.g.~\cite{Xu98}): 
The quantum vanishing ideal $\qideal_C$ is a left ideal corresponding
to the \emph{coisotropic} sub-manifold $C = J^{-1} (\{0\})$. This would
also support the `coisotropic creed' formulated by Lu in \cite{Lu93}.
Furthermore, one may speculate whether a connected \emph{Lagrangian}
sub-manifold corresponds to a left ideal $\qideal$ whose Lie idealiser 
is equal to $\qideal \oplus \mathbb C[[\lambda]] 1$. This seems to be
reasonable in view of the situation for cotangent bundles:
As discussed in detail in \cite{BNW98a,BNW99a,BNPW98} for any
cotangent bundle $T^*Q$ there is a homogeneous star product of Weyl
type such that the integration over the configuration space with
respect to a positive density is a formally positive functional whose
Gel'fand ideal (which is a left ideal) characterises the configuration 
space as embedded Lagrangian sub-manifold. Here one can rather easily
show that the Lie idealiser of this Gel'fand ideal can be obtained by
adding the multiples of the function $1$ if the configuration space
$Q$ is connected. We shall not go into details but mention that the
above characterisation is still true for projectable Lagrangian
sub-manifolds $L$ if one takes the Gel'fand ideal of  
a formally positive functional having support on $L$ which also
induces the WKB expansion.

%
%

\section* {Acknowledgements}

We would like to thank Peter Gl{\"o}{\ss}ner, 
Joachim Schirmer, and Gijs Tuynman for many clarifying discussions. Moreover,
we are very grateful to G.Halbout, M. Semenov-Tyan-Shanskii,
J.~Stasheff, A.~Weinstein, and to the referee for many useful
discussions, suggestions, and constructive criticisms concerning the
first version of this article.

\end{document}